\documentclass[a4paper,12pt]{article} 
\usepackage{amsmath} 
\usepackage{amsopn} 
\usepackage{amssymb} 
\usepackage[mathscr]{eucal} 
\usepackage{theorem} 
\usepackage{enumerate} 

\theorembodyfont{\itshape} 
\theoremstyle{plain}
  \newtheorem{thm}{Theorem}[section] 
   
  \newtheorem{lem}[thm]{Lemma}

\theorembodyfont{\rmfamily} 
\theoremstyle{plain}
  \newtheorem{rem}{Remark}[section] 
  \newtheorem{ex}{Example}[section] 

\renewcommand{\theequation}%
           {\thesection.\arabic{equation}}

\begin{document} 

\begin{center} 
{\LARGE Surfaces with flat normal connection} 

\vspace{1mm} 

{\LARGE in 4-dimensional space forms} 

\vspace{6mm} 

{\Large Naoya {\sc Ando} and Ryusei {\sc Hatanaka}} 
\end{center} 

\vspace{3mm} 

\begin{quote} 
{\footnotesize \it Abstract} \ 
{\footnotesize 
Let $N$ be a Riemannian, Lorentzian or neutral $4$-dimensional space form 
with constant sectional curvature $L_0$. 
In this paper, 
noticing the linearly dependent condition, 
we obtain characterizations of space-like surfaces in $N$ 
with flat normal connection and parallel normal vector fields. 
In addition, 
we obtain a generic characterization of space-like surfaces in $N$ 
with flat normal connection and $K\equiv L_0$ 
which do not admit any parallel normal vector fields. 
For time-like surfaces in $N$ with flat normal connection, 
we obtain analogous results.}
\end{quote} 

\vspace{3mm} 

\section{Introduction}\label{sect:intro} 

\setcounter{equation}{0} 

Let $N$ be a $4$-dimensional Riemannian space form  
with constant sectional curvature $L_0$. 
Let $M$ be a Riemann surface 
and $F:M\longrightarrow N$ a conformal immersion of $M$ into $N$. 
The normal connection $\nabla^{\perp}$ of $F$ is said to be \textit{flat\/} 
if the curvature tensor $R^{\perp}$ of $\nabla^{\perp}$ vanishes. 
See \cite{chen}, \cite{hoffman} for Hoffman surfaces, 
which are surfaces with flat normal connection 
and parallel mean curvature vector. 
Based on \cite{thas}, it is observed in \cite{enomoto} 
that a flat surface in Euclidean $4$-space $E^4$ 
with flat normal connection and nondegenerate Gauss map is represented as 
the Riemannian product of two plane curves. 
In particular, there exist flat surfaces in $E^4$ 
with flat normal connection and no parallel normal vector fields. 
See \cite{DT}, \cite{GM}, \cite{MS} 
for constructions of flat surfaces in $E^4$ 
with flat normal connection. 

In the present paper, 
noticing the \textit{linearly dependent condition\/} 
(see Subsection~\ref{subsect:Rpnvf}), 
which is related to existence of parallel normal vector fields, 
we will study space-like or time-like surfaces 
in Riemannian, Lorentzian or neutral $4$-dimensional space forms 
with flat normal connection. 

Let $N$, $M$, $F$ be as in the beginning of this section. 
If $F$ has a parallel normal vector field 
with respect to the Levi-Civita connection of the metric of $N$, 
that is, 
if the image of $F$ is contained in a $3$-dimensional totally geodesic 
hypersurface of $N$, 
then the normal connection $\nabla^{\perp}$ is flat and   
the second fundamental form $\sigma$ of $F$ satisfies 
the linearly dependent condition. 
Let $K$ be the curvature of the induced metric by $F$. 
We will see that 
if $F$ is a conformal immersion 
with $R^{\perp} \equiv 0$ and $K\not= L_0$ 
such that $\sigma$ satisfies the linearly dependent condition, 
then $F$ has a parallel normal vector field  
((a) of Theorem~\ref{thm:pnvf}). 
Therefore, for a conformal immersion $F$ 
with $R^{\perp} \equiv 0$ and $K\not= L_0$, 
the existence of a parallel normal vector field is equivalent to 
the linearly dependent condition of $\sigma$. 

For a conformal immersion $F$ 
with $R^{\perp} \equiv 0$ and $K\equiv L_0$ 
such that $\sigma$ satisfies the linearly dependent condition, 
we will obtain a criterion 
for the existence of a parallel normal vector field 
((b) of Theorem~\ref{thm:pnvf}), and 
we will observe that there exist conformal immersions 
with $R^{\perp} \equiv 0$, $K\equiv L_0$ 
and no parallel normal vector fields 
such that the second fundamental forms satisfy 
the linearly dependent condition (Theorem~\ref{thm:ldnotpnvf}). 
Therefore, for a conformal immersion $F$ 
with $R^{\perp} \equiv 0$ and $K\equiv L_0$, 
the linearly dependent condition of $\sigma$ does not necessarily mean 
the existence of a parallel normal vector field. 

There exist conformal immersions of $M$ into $N$ 
with $R^{\perp} \equiv 0$ and $K\equiv L_0$ 
such that the second fundamental forms do not satisfy  
the linearly dependent condition. 
Such immersions have no parallel normal vector fields. 
By the equations of Gauss and Ricci, 
the two conditions $R^{\perp} \equiv 0$ and $K\equiv L_0$ just mean 
the relations among the coefficients of $\sigma$ given in \eqref{WXYZ0}, 
and in \cite{ando10}, 
it is observed that \eqref{WXYZ0} is equivalent to 
the degeneracy of the twistor lifts of $F$. 
Based on \eqref{WXYZ0} and Lemma~\ref{lem:XYXY}, 
and studying the equations of Codazzi, 
we will obtain a generic characterization of 
a conformal immersion $F$ with $R^{\perp} \equiv 0$ and $K\equiv L_0$ 
such that $\sigma$ does not satisfy 
the linearly dependent condition (Theorem~\ref{thm:notldL0}). 
In addition, in the case of $L_0 =0$, 
we will obtain a more concrete characterization 
by an over-determined system of polynomial type with degree at most two 
(Theorem~\ref{thm:notld}). 
See \cite{ando7} for over-determined systems of polynomial type. 

Let $N$ be a $4$-dimensional neutral space form  
with constant sectional curvature $L_0$. 
Let $M$ be a Riemann surface 
and $F:M\longrightarrow N$ a space-like and conformal immersion 
of $M$ into $N$. 
Then noticing the equations of Gauss, Codazzi and Ricci, 
we can obtain the same results as in the case 
where $N$ is Riemannian (see Section~\ref{sect:sn}). 

Let $N$ be as in the previous paragraph. 
Let $M$ be a Lorentz surface 
and $F:M\longrightarrow N$ a time-like and conformal immersion 
of $M$ into $N$. 
Suppose that $F$ has a parallel normal vector field 
with respect to the Levi-Civita connection of the neutral metric of $N$. 
Then this vector field is identically space-like, time-like or 
light-like. 
If it is space-like, 
then the image of $F$ is contained in 
a $3$-dimensional totally geodesic hypersurface of $N$ 
which has a nondegenerate and indefinite metric 
with signature $(1, 2)$; 
if it is time-like, 
then the image of $F$ is contained in 
a $3$-dimensional totally geodesic Lorentzian hypersurface of $N$; 
if it is light-like, 
then the curvature $K$ of the induced metric by $F$ is identically equal to 
$L_0$ by the equation of Gauss (see \eqref{gausstn}), 
and the image of $F$ is contained in 
a $3$-dimensional totally geodesic hypersurface of $N$ 
with degenerate metric. 
For time-like and conformal immersions of $M$ into $N$ 
with space-like or time-like parallel normal vector fields, 
we will obtain an analogous result to Theorem~\ref{thm:pnvf} 
(see Theorem~\ref{thm:pnvftn}), 
and referring to Subsection~\ref{subsect:Rpnvf}, 
we can obtain analogous results to Theorem~\ref{thm:ldnotpnvf}. 
We will characterize time-like and conformal immersions of $M$ into $N$ 
with light-like parallel normal vector fields 
(Theorem~\ref{thm:pnvftnll}). 
The conformal Gauss maps of time-like surfaces in 
$3$-dimensional Lorentzian space forms of Willmore type 
with vanishing paraholomorphic quartic differential have 
light-like parallel normal vector fields (\cite{ando4}). 
See Remark~\ref{rem:K=L_0} in Subsection~\ref{subsect:tnpnvf} 
for characterizations of 
a time-like and conformal immersion $F:M\longrightarrow N$ 
with zero mean curvature vector and $K\equiv L_0$. 
We will obtain analogous results to 
Theorems~\ref{thm:notldL0} and \ref{thm:notld} 
for time-like and conformal immersions of $M$ into $N$ 
with $R^{\perp} \equiv 0$, $K\equiv L_0$ 
and no parallel normal vector fields 
(Theorems~\ref{thm:notldL0tn}, \ref{thm:notldtn}). 
The two conditions $R^{\perp} \equiv 0$ and $K\equiv L_0$ just mean 
the relations in \eqref{WXYZ0tn} analogous to \eqref{WXYZ0}, 
and \eqref{WXYZ0tn} is equivalent to 
the degeneracy of the time-like twistor lifts of $F$ (\cite{ando10}). 

Let $N$ be a $4$-dimensional Lorentzian space form  
with constant sectional curvature $L_0$. 
Let $M$ be a Riemann surface 
and $F:M\longrightarrow N$ a space-like and conformal immersion 
of $M$ into $N$. 
If $F$ has a parallel normal vector field 
with respect to the Levi-Civita connection of the Lorentzian metric of $N$, 
then as in the previous paragraph, 
this vector field is identically space-like, time-like or 
light-like. 
For space-like and conformal immersions of $M$ into $N$ 
with space-like or time-like parallel normal vector fields, 
we can obtain analogous results to Theorems~\ref{thm:ldnotpnvf} 
and \ref{thm:pnvftn} for this situation. 
In the case where the immersions have 
light-like parallel normal vector fields, 
we can also obtain an analogous result to Theorem~\ref{thm:pnvftnll}. 
The conformal Gauss maps of Willmore surfaces in 
$3$-dimensional Riemannian space forms  
with vanishing holomorphic quartic differential have 
light-like parallel normal vector fields (\cite{bryant2}, \cite{ando3}). 
See Remark~\ref{rem:K=L_0sL} in Subsection~\ref{subsect:sLGCR} 
for characterizations of 
a space-like and conformal immersion $F:M\longrightarrow N$ 
with zero mean curvature vector and $K\equiv L_0$. 
We will obtain analogous results to 
Theorems~\ref{thm:notldL0} and \ref{thm:notld} 
for space-like and conformal immersions of $M$ into $N$ 
with $K\equiv L_0$, $R^{\perp} \equiv 0$ 
and no parallel normal vector fields 
(Theorems~\ref{thm:notldL0sL}, \ref{thm:notldsL}), 
expressing the equations of Gauss, Codazzi and Ricci 
as in \eqref{gaussriccisL}, \eqref{codazziWXYZsL} 
by complex-valued functions related to the second fundamental forms. 
The two conditions $R^{\perp} \equiv 0$ and $K\equiv L_0$ just mean 
the relation \eqref{WXYZ0sL} analogous to \eqref{WXYZ0}, 
and \eqref{WXYZ0sL} is equivalent to 
the degeneracy of the complex twistor lifts of $F$ (\cite{ando10}). 

Let $N$ be as in the previous paragraph. 
Let $M$ be a Lorentz surface 
and $F:M\longrightarrow N$ a time-like and conformal immersion 
of $M$ into $N$. 
If $F$ has a parallel normal vector field 
with respect to the Levi-Civita connection, 
then this vector field is identically space-like. 
For time-like and conformal immersions of $M$ into $N$ 
with parallel normal vector fields, 
we can obtain analogous results to Theorems~\ref{thm:pnvf} 
and \ref{thm:ldnotpnvf} for this situation. 
We will obtain analogous results to 
Theorems~\ref{thm:notldL0sL} and \ref{thm:notldsL} 
for time-like and conformal immersions of $M$ into $N$ 
with $R^{\perp} \equiv 0$, $K\equiv L_0$ 
and no parallel normal vector fields 
(Theorems~\ref{thm:notldL0tL}, \ref{thm:notldtL}). 
The two conditions $R^{\perp} \equiv 0$ and $K\equiv L_0$ just mean 
the relation \eqref{WXYZ0tL} analogous to \eqref{WXYZ0sL}, 
and \eqref{WXYZ0tL} is equivalent to 
the degeneracy of the complex twistor lifts of $F$ (\cite{ando10}). 

\section{Surfaces with flat normal connection 
in Riemannian 4-dimensional space forms}\label{sect:R} 

\setcounter{equation}{0} 

\subsection{The Gauss-Codazzi-Ricci equations}\label{subsect:RGCR} 

Let $N$ be a 4-dimensional Riemannian space form 
with constant sectional curvature $L_0$. 
If $L_0 =0$, 
then we can suppose $N=E^4$; 
if $L_0 >0$, 
then we can suppose $N=\{ x\in E^5 \ | \ \langle x, x\rangle =1/L_0 \}$; 
if $L_0 <0$, 
then we can suppose $N=\{ x\in E^5_1 \ | \ \langle x, x\rangle =1/L_0 \}$, 
where $\langle \ , \ \rangle$ is the metric of $E^5$ or $E^5_1$.  
Let $h$ be the metric of $N$. 
Let $M$ be a Riemann surface 
and $F:M\longrightarrow N$ a conformal immersion of $M$ into $N$. 
Let $(u, v)$ be local isothermal coordinates of $M$ 
compatible with the complex structure of $M$. 
Then the induced metric $g$ on $M$ by $F$ is represented 
as $g=e^{2\lambda} (du^2 +dv^2 )$ 
for a real-valued function $\lambda$. 
Then we have 
$$h(dF(\partial_u ), dF(\partial_u )) 
 =h(dF(\partial_v ), dF(\partial_v )), \quad  
  h(dF(\partial_u ), dF(\partial_v ))=0$$ 
($\partial_u :=\partial /\partial u$, 
 $\partial_v :=\partial /\partial v$). 
Let $\Tilde{\nabla}$ denote the Levi-Civita connection 
of $E^4$, $E^5$ or $E^5_1$ according to $L_0 =0$, $>0$ or $<0$. 
We set $T_1 :=dF(\partial_u )$, $T_2 :=dF(\partial_v )$. 
Let $N_1$, $N_2$ be normal vector fields of $F$ satisfying 
$$h(N_1 , N_1 )=h(N_2 , N_2 )=e^{2\lambda} , \quad 
  h(N_1 , N_2 )=0.$$ 
Then we have 
\begin{equation} 
\begin{split} 
\Tilde{\nabla}_{T_1} (T_1 \ T_2 \ N_1 \ N_2 \ F) 
& =                  (T_1 \ T_2 \ N_1 \ N_2 \ F)S, \\ 
\Tilde{\nabla}_{T_2} (T_1 \ T_2 \ N_1 \ N_2 \ F) 
& =                  (T_1 \ T_2 \ N_1 \ N_2 \ F)T, 
\end{split} 
\label{dt1t2} 
\end{equation} 
where 
\begin{equation} 
\begin{split} 
S & =\left[ 
     \begin{array}{ccccc} 
      \lambda_u         & \lambda_v & -\alpha_1  & -\beta_1   & 1 \\ 
     -\lambda_v         & \lambda_u & -\alpha_2  & -\beta_2   & 0 \\ 
      \alpha_1          & \alpha_2  &  \lambda_u & -\mu_1     & 0 \\ 
      \beta_1           & \beta_2   &  \mu_1     &  \lambda_u & 0 \\ 
      -L_0 e^{2\lambda} &  0        &   0        &   0        & 0 
       \end{array} 
     \right] , \\ 
T & =\left[ 
     \begin{array}{ccccc} 
      \lambda_v & -\lambda_u        & -\alpha_2  & -\beta_2   & 0 \\ 
      \lambda_u &  \lambda_v        & -\alpha_3  & -\beta_3   & 1 \\ 
      \alpha_2  &  \alpha_3         &  \lambda_v & -\mu_2     & 0 \\ 
      \beta_2   &  \beta_3          &  \mu_2     &  \lambda_v & 0 \\ 
       0        & -L_0 e^{2\lambda} &   0        &   0        & 0 
       \end{array} 
    \right] , 
\end{split} 
\label{ST} 
\end{equation} 
and $\alpha_k$, $\beta_k$ ($k=1, 2, 3$) and $\mu_l$ ($l=1, 2$) are 
real-valued functions. 
From \eqref{dt1t2}, we obtain $S_v -T_u =ST-TS$. 
This is equivalent to 
the system of the equations of Gauss, Codazzi and Ricci. 
The equation of Gauss is given by 
\begin{equation} 
  \lambda_{uu} +\lambda_{vv} +L_0 e^{2\lambda} 
=-\alpha_1 \alpha_3 -\beta_1 \beta_3 +\alpha^2_2 +\beta^2_2 . 
\label{gauss}
\end{equation}
The equations of Codazzi are given by 
\begin{equation} 
\begin{split}
  (\alpha_1 )_v -(\alpha_2 )_u 
& =\alpha_2 \lambda_u +\alpha_3 \lambda_v 
   -\beta_2 \mu_1     +\beta_1  \mu_2  , \\ 
   (\alpha_2 )_v -(\alpha_3 )_u 
& =-\alpha_1 \lambda_u -\alpha_2 \lambda_v 
   -\beta_3  \mu_1     +\beta_2  \mu_2  , \\ 
   (\beta_1 )_v -(\beta_2 )_u 
& = \beta_2  \lambda_u +\beta_3  \lambda_v 
   +\alpha_2 \mu_1     -\alpha_1 \mu_2  , \\ 
   (\beta_2 )_v -(\beta_3 )_u 
& =-\beta_1  \lambda_u -\beta_2  \lambda_v 
   +\alpha_3 \mu_1     -\alpha_2 \mu_2 .  
\end{split} 
\label{codazzi} 
\end{equation} 
The equation of Ricci is given by 
\begin{equation} 
  (\mu_1 )_v -(\mu_2 )_u = \alpha_1 \beta_2 -\alpha_2 \beta_1 
                          +\alpha_2 \beta_3 -\alpha_3 \beta_2 . 
\label{ricci}
\end{equation}
Let $\nabla$ be the Levi-Civita connection of $h$. 
Then for the normal connection $\nabla^{\perp}$ of $F$ in $N$, 
$\nabla^{\perp} N_l$ ($l\in \{ 1, 2\}$) is given by 
the normal component of $\nabla N_l$. 
Let $R^{\perp}$ be the curvature tensor of $\nabla^{\perp}$. 
Then the condition that $\nabla^{\perp}$ is flat just means 
$R^{\perp} \equiv 0$, and 
$R^{\perp} \equiv 0$ is equivalent to $(\mu_2 )_u =(\mu_1 )_v$. 

\subsection{Surfaces with \mbox{\boldmath{$R^{\perp} \equiv 0$}} and 
parallel normal vector fields}\label{subsect:Rpnvf} 

Suppose that $\nabla^{\perp}$ is flat. 
Then there exists a function $\gamma$ satisfying 
$\gamma_u =\mu_1$, $\gamma_v =\mu_2$. 
In addition, we obtain 
\begin{equation} 
 \alpha_1 \beta_2 -\alpha_2 \beta_1 
+\alpha_2 \beta_3 -\alpha_3 \beta_2 =0. 
\label{riccifnc} 
\end{equation} 
We say that the second fundamental form $\sigma$ of $F$ satisfies 
the \textit{linearly dependent condition\/} 
if $F$ satisfies 
\begin{equation} 
 \cos \theta (\alpha_1 , \alpha_2 , \alpha_3 ) 
+\sin \theta (\beta_1  , \beta_2  , \beta_3  ) 
=0
\label{f} 
\end{equation} 
for a function $\theta$ of $u$, $v$. 

We will prove 

\begin{thm}\label{thm:pnvf} 
Suppose that $\nabla^{\perp}$ is flat. 
\begin{itemize}
\item[{\rm (a)}]{Suppose that the curvature $K$ of $g$ is nowhere equal 
to $L_0$. 
Then on a neighborhood of each point of $M$, 
$F$ has a parallel normal vector field 
if and only if the second fundamental form $\sigma$ of $F$ satisfies 
the linearly dependent condition. 
In addition, if one of these conditions holds, 
then $\gamma +\theta$ is constant.} 
\item[{\rm (b)}]{Suppose $K\equiv L_0$. 
Then on a neighborhood of each point of $M$, 
$F$ has a parallel normal vector field 
if and only if $\sigma$ satisfies the linearly dependent condition 
so that $\gamma +\theta$ is constant.} 
\end{itemize} 
If $F$ has a parallel normal vector field, 
then it is given by 
$e^{-\lambda} (\cos \theta N_1 +\sin \theta N_2 )$. 
\end{thm} 

\vspace{3mm} 

\par\noindent 
\textit{Proof} \ 
If on a neighborhood of each point of $M$, 
$F$ has a parallel normal vector field, 
then we can suppose that $e^{-\lambda} N_1$ is parallel, 
and then $\gamma$ is constant and $\theta =0$ satisfies \eqref{f}. \\ 
Suppose $K\not= L_0$. 
Suppose that there exists a function $\theta$ satisfying \eqref{f} 
and $\sin \theta \not= 0$. 
Applying \eqref{f} to \eqref{codazzi}, we obtain 
\begin{equation} 
\begin{split} 
    \alpha_2 \gamma_u -\alpha_1 \gamma_v 
& =-\alpha_2 \theta_u +\alpha_1 \theta_v \\ 
    \alpha_3 \gamma_u -\alpha_2 \gamma_v 
& =-\alpha_3 \theta_u +\alpha_2 \theta_v . 
\end{split} 
\label{codazzifnc0} 
\end{equation} 
We can rewrite \eqref{codazzifnc0} into 
\begin{equation} 
\left[ \begin{array}{cc} 
        \alpha_1 & \alpha_2 \\ 
        \alpha_2 & \alpha_3  
         \end{array} 
\right] 
\left[ \begin{array}{c} 
       -(\gamma +\theta )_v \\ 
        (\gamma +\theta )_u 
         \end{array} 
\right] 
=
\left[ \begin{array}{c} 
        0 \\ 
        0 
         \end{array} 
\right] . 
\label{codazzifnc} 
\end{equation} 
By \eqref{gauss}, \eqref{f} and $K\not= L_0$, 
we obtain $\alpha_1 \alpha_3 -\alpha^2_2 \not= 0$. 
Therefore the matrix in \eqref{codazzifnc} is nonsingular 
and $\gamma +\theta$ is constant. 
Then we obtain 
\begin{equation} 
           \nabla (\cos \theta N_1 +\sin \theta N_2 ) 
=d\lambda \otimes (\cos \theta N_1 +\sin \theta N_2 ). 
\label{xi} 
\end{equation} 
This means that $e^{-\lambda} (\cos \theta N_1 +\sin \theta N_2 )$ is 
parallel. 
Hence we obtain (a) of Theorem~\ref{thm:pnvf}. 
\\ 
Suppose $K\equiv L_0$. 
Suppose that there exists a function $\theta$ satisfying 
\eqref{f} and that $\gamma +\theta$ is constant. 
Then we obtain \eqref{xi} 
and $e^{-\lambda} (\cos \theta N_1 +\sin \theta N_2 )$ is parallel. 
Hence we obtain (b) of Theorem~\ref{thm:pnvf}. 
\hfill 
$\square$ 

\begin{rem} 
Suppose one of the conditions in (b) of Theorem~\ref{thm:pnvf}. 
Since $K\equiv L_0$, 
we obtain $\alpha_1 \alpha_3 -\alpha^2_2 =0$ 
or        $\beta_1  \beta_3  -\beta^2_2  =0$. 
We suppose 
\begin{equation*} 
\alpha_1 \alpha_3 -\alpha^2_2 =0, \quad 
(\alpha_1 , \alpha_2 , \alpha_3 )\not= (0, 0, 0). 
\end{equation*} 
Suppose $\alpha_2 \equiv 0$. 
Then $\alpha_1 \equiv 0$ or $\alpha_3 \equiv 0$. 
Suppose $\alpha_3 \equiv 0$. 
Then $\alpha_1 \not= 0$, and then $\sin \theta \not= 0$. 
We can suppose $\alpha_1 \sin \theta >0$. 
From \eqref{codazzi}, 
we obtain $\lambda_u \equiv 0$ 
and $\alpha_1 /\sin \theta$ is a function of one variable $u$. 
In the following, suppose $\alpha_2 \not= 0$. 
Then $\alpha_1 \not= 0$ and $\alpha_3 \not= 0$. 
We suppose $\alpha_1 >0$ and $\sin \theta >0$. 
We set $\rho :=\alpha_2 /\alpha_1 =\alpha_3 /\alpha_2$. 
Applying \eqref{f} to \eqref{codazzi}, we obtain \eqref{codazzifnc0}. 
Applying \eqref{codazzifnc0} to the first two relations 
in \eqref{codazzi}, we obtain 
\begin{equation} 
\begin{split} 
      \left( \dfrac{\alpha_1}{\sin \theta} \right)_v 
     -\left( \dfrac{\alpha_2}{\sin \theta} \right)_u 
  & = \dfrac{1}{\sin \theta} (\alpha_2 \lambda_u +\alpha_3 \lambda_v ), \\ 
      \left( \dfrac{\alpha_2}{\sin \theta} \right)_v 
     -\left( \dfrac{\alpha_3}{\sin \theta} \right)_u 
  & =-\dfrac{1}{\sin \theta} (\alpha_1 \lambda_u +\alpha_2 \lambda_v ). 
\end{split} 
\label{codazzifnc2} 
\end{equation} 
By \eqref{codazzifnc2} and $\rho =\alpha_2 /\alpha_1 =\alpha_3 /\alpha_2$, 
we obtain 
\begin{equation*}  
   \left( \dfrac{\alpha_1}{\sin \theta} \sqrt{\rho^2 +1} \right)_v 
  =\left( \dfrac{\alpha_2}{\sin \theta} \sqrt{\rho^2 +1} \right)_u . 
\end{equation*} 
Therefore there exists a function $\phi$ of two variables $u$, $v$ 
satisfying $\phi^2_u +\phi^2_v \not= 0$ and 
\begin{equation} 
  \phi_u =\dfrac{\alpha_1}{\sin \theta} \sqrt{\rho^2 +1} , \quad 
  \phi_v =\dfrac{\alpha_2}{\sin \theta} \sqrt{\rho^2 +1} . 
\label{phiuphiv} 
\end{equation} 
From \eqref{phiuphiv}, we obtain $\phi_v =\rho \phi_u$. 
Applying this to the first relation in \eqref{phiuphiv} 
and noticing $\alpha_1 >0$, we obtain 
\begin{equation*} 
  \dfrac{\alpha_1}{\sin \theta} =\dfrac{\phi^2_u}{\sqrt{\phi^2_u +\phi^2_v}} . 
\label{alphaphi} 
\end{equation*} 
Therefore we obtain 
\begin{equation} 
(\alpha_1 , \alpha_2 , \alpha_3 ) 
=\dfrac{\sin \theta}{\sqrt{\phi^2_u +\phi^2_v}} 
(\phi^2_u , \phi_u \phi_v , \phi^2_v ).  
\label{alpha} 
\end{equation} 
Applying \eqref{alpha} to one of the relations in \eqref{codazzifnc2}, 
we obtain 
\begin{equation} 
    \phi^2_v \phi_{uu} -2\phi_u \phi_v \phi_{uv} +\phi^2_u \phi_{vv} 
  +(\phi^2_u +\phi^2_v )(\phi_u \lambda_u +\phi_v \lambda_v )
  =0. 
\label{phi} 
\end{equation}
\end{rem} 

\subsection{Surfaces with \mbox{\boldmath{$R^{\perp} \equiv 0$}} and 
                          \mbox{\boldmath{$K\equiv L_0$}} 
which do not admit any parallel normal vector fields}\label{subsect:Rnpnvf} 

Let $\lambda$ be a function satisfying $\lambda_{uu} +\lambda_{vv} +L_0 e^{2\lambda} =0$. 
Let $\phi$ be a function satisfying \eqref{phi}. 
Let $\theta$ be a function satisfying $\sin \theta \cos \theta \not= 0$. 
Let $\alpha_k$, $\beta_k$ ($k=1, 2, 3$) be functions 
given by \eqref{f} and \eqref{alpha}. 
We set $\gamma :=-\theta +\xi (\phi)$, 
where $\xi$ is a function of one variable. 
Then we can choose $\xi$ so that $\gamma +\theta$ is not constant, and 
$\lambda$, $\alpha_k$, $\beta_k$ ($k=1, 2, 3$), 
$\mu_1 :=\gamma_u$, $\mu_2 :=\gamma_v$ satisfy 
\eqref{gauss}, \eqref{codazzi}, \eqref{ricci}. 
Therefore we obtain 

\begin{thm}\label{thm:ldnotpnvf} 
There exist surfaces in $N$ with $K\equiv L_0$ and $R^{\perp} \equiv 0$ 
satisfying 
\begin{itemize}
\item[{\rm (a)}]{the second fundamental forms satisfy 
the linearly dependent condition\/$;$} 
\item[{\rm (b)}]{there do not exist any parallel normal vector fields.}  
\end{itemize} 
\end{thm} 

In the following, 
we study a conformal immersion $F:M\longrightarrow N$ of $M$ into $N$ 
with $R^{\perp} \equiv 0$ and $K\equiv L_0$ 
such that the second fundamental form of $F$ does not satisfy 
the linearly dependent condition. 

We set 
\begin{equation} 
\begin{array}{rcr} 
W_{\pm} :=\alpha_2 \pm \beta_1  , & \ & X_{\pm} :=\alpha_2 \pm \beta_3 , \\ 
Y_{\pm} :=\beta_2  \pm \alpha_1 , & \ & Z_{\pm} :=\beta_2  \pm \alpha_3 . 
\end{array} 
\label{WXYZ} 
\end{equation} 
Then we have 
\begin{equation} 
W_+ +W_- =X_+ +X_- , \quad 
Y_+ +Y_- =Z_+ +Z_- . 
\label{W=XY=Z} 
\end{equation} 
By \eqref{gauss} and \eqref{ricci}, we obtain 
\begin{equation} 
 W_{\mp} X_{\pm} +Y_{\pm} Z_{\mp} 
=\lambda_{uu} +\lambda_{vv} +L_0 e^{2\lambda} 
 \pm ((\mu_1 )_v -(\mu_2 )_u ). 
\label{gaussricci} 
\end{equation} 
We can rewrite \eqref{codazzi} into 
\begin{equation} 
\begin{split}
   (Y_+ )_v -(X_+ )_u 
& = W_- (\lambda_u -\mu_2 )-Z_- (\lambda_v +\mu_1 ), \\ 
   (Y_- )_v +(X_- )_u 
& =-W_+ (\lambda_u +\mu_2 )-Z_+ (\lambda_v -\mu_1 ), \\ 
   (W_+ )_v -(Z_+ )_u 
& = Y_- (\lambda_u +\mu_2 )-X_- (\lambda_v -\mu_1 ), \\ 
   (W_- )_v +(Z_- )_u 
& =-Y_+ (\lambda_u -\mu_2 )-X_+ (\lambda_v +\mu_1 ).  
\end{split} 
\label{codazziWXYZ} 
\end{equation} 

Suppose $R^{\perp} \equiv 0$ and $K\equiv L_0$. 
Then from \eqref{gaussricci}, 
we obtain 
\begin{equation} 
W_{\mp} X_{\pm} +Y_{\pm} Z_{\mp} =0. 
\label{WXYZ0} 
\end{equation} 
We suppose that there exist functions $k_{\pm}$ 
satisfying 
\begin{equation} 
(W_- , Z_- )=k_+ (-Y_+ , X_+ ), \quad 
(W_+ , Z_+ )=k_- (-Y_- , X_- ). 
\label{kpm} 
\end{equation} 
Then $k_{\pm} \equiv 0$ if and only if 
$\beta_1$, $\beta_2$, $\alpha_2$, $\alpha_3$ vanish. 
The product of two plane curves in $E^4$ satisfies this condition 
and such surfaces are studied in \cite{enomoto}. 
In the following, suppose $k_{\pm} \not= 0$. 
We will prove 

\begin{lem}\label{lem:fXY} 
There exist functions $f_+$, $f_-$ satisfying 
\begin{equation} 
\begin{array}{lcl} 
 X_+ = \dfrac{(f_+ )_v}{\sqrt{1+k^2_+}} , & \ &  
 Y_+ = \dfrac{(f_+ )_u}{\sqrt{1+k^2_+}} ,   \\ 
 X_- =-\dfrac{(f_- )_v}{\sqrt{1+k^2_-}} , & \ &  
 Y_- = \dfrac{(f_- )_u}{\sqrt{1+k^2_-}} . 
  \end{array} 
\label{fpm}
\end{equation} 
\end{lem} 

\vspace{3mm} 

\par\noindent 
\textit{Proof} \ 
Applying \eqref{kpm} to \eqref{codazziWXYZ}, we obtain 
\begin{equation} 
\begin{split}
   (Y_+ )_v -(X_+ )_u 
& =-k_+ Y_+ (\lambda_u -\gamma_v ) 
   -k_+ X_+ (\lambda_v +\gamma_u ), \\ 
  -(Y_- )_v -(X_- )_u 
& =-k_- Y_- (\lambda_u +\gamma_v ) 
   +k_- X_- (\lambda_v -\gamma_u ), \\ 
  -(k_- Y_- )_v -(k_- X_- )_u 
& = Y_- (\lambda_u +\gamma_v )-X_- (\lambda_v -\gamma_u ), \\ 
  -(k_+ Y_+ )_v +(k_+ X_+ )_u 
& =-Y_+ (\lambda_u -\gamma_v )-X_+ (\lambda_v +\gamma_u ).  
\end{split} 
\label{codazzikXY} 
\end{equation} 
From the first and fourth equations in \eqref{codazzikXY}, 
we obtain 
\begin{equation} 
(k_+ )_u X_+ -(k_+ )_v Y_+ 
=-(1+k^2_+ )(X_+ (\lambda_v +\gamma_u )+Y_+ (\lambda_u -\gamma_v ));  
\label{kXY+} 
\end{equation} 
from the second and third equations in \eqref{codazzikXY}, 
we obtain 
\begin{equation} 
(k_- )_u X_- +(k_- )_v Y_- 
= (1+k^2_- )(X_- (\lambda_v -\gamma_u )-Y_- (\lambda_u +\gamma_v )). 
\label{kXY-} 
\end{equation} 
From the first equation in \eqref{codazzikXY} and \eqref{kXY+}, 
we obtain 
\begin{equation*} 
 \left( \sqrt{1+k^2_+} X_+ \right)_u 
=\left( \sqrt{1+k^2_+} Y_+ \right)_v ; 
\label{kXY+2} 
\end{equation*} 
from the second equation in \eqref{codazzikXY} and \eqref{kXY-}, 
we obtain 
\begin{equation*} 
  \left( \sqrt{1+k^2_-} X_- \right)_u 
=-\left( \sqrt{1+k^2_-} Y_- \right)_v . 
\label{kXY-2} 
\end{equation*} 
Therefore there exist functions $f_{\pm}$ satisfying \eqref{fpm}. 
Hence we obtain Lemma~\ref{lem:fXY}. 
\hfill 
$\square$ 

\vspace{3mm} 

We will prove 

\begin{lem}\label{lem:f+2f-2} 
The functions $f_+$, $f_-$ satisfy 
\begin{equation} 
(f_+ )^2_u + (f_+ )^2_v =(f_- )^2_u + (f_- )^2_v . 
\label{grad} 
\end{equation} 
\end{lem} 

\vspace{3mm} 

\par\noindent 
\textit{Proof} \ 
Since \eqref{riccifnc} is equivalent to 
\begin{equation}  
 W^2_+ +X^2_- +Y^2_- +Z^2_+ 
=W^2_- +X^2_+ +Y^2_+ +Z^2_- , 
\label{X2Y2Z2W2} 
\end{equation} 
from \eqref{kpm}, \eqref{fpm} and \eqref{X2Y2Z2W2}, 
we obtain \eqref{grad}. 
\hfill 
$\square$ 

\vspace{3mm} 
 
We will prove 

\begin{lem}\label{lem:XYXY} 
If $X^2_+ Y^2_- -X^2_- Y^2_+ \not= 0$, 
then the second fundamental form of $F$ does not satisfy 
the linearly dependent condition. 
\end{lem} 

\vspace{3mm} 

\par\noindent 
\textit{Proof} \ 
Suppose that a function $\theta$ satisfies \eqref{f}. 
Then $\alpha_1 \beta_2 -\alpha_2 \beta_1 =0$ 
and  $\alpha_2 \beta_3 -\alpha_3 \beta_2 =0$. 
By \eqref{kpm}, we obtain 
\begin{equation*} 
\begin{split} 
    4(\alpha_1 \beta_2 -\alpha_2 \beta_1 ) 
& =(1+k^2_+ )Y^2_+ -(1+k^2_- )Y^2_- , \\ 
    4(\alpha_2 \beta_3 -\alpha_3 \beta_2 ) 
& =(1+k^2_+ )X^2_+ -(1+k^2_- )X^2_- . 
\end{split} 
\end{equation*} 
Therefore we obtain $X^2_+ Y^2_- -X^2_- Y^2_+ =0$. 
Hence we obtain Lemma~\ref{lem:XYXY}. 
\hfill 
$\square$ 

\vspace{3mm} 

In the following, we suppose 
\begin{equation} 
  X^2_+ Y^2_- -X^2_- Y^2_+ 
=(X_+ Y_- +X_- Y_+ )(X_+ Y_- -X_- Y_+ )\not= 0. 
\label{X+Y-X-Y+} 
\end{equation} 
We set 
\begin{equation} 
\begin{split} 
A & :=(f_+ )_v (f_- )_u +(f_- )_v (f_+ )_u , \quad 
A'  :=(f_+ )_v (f_- )_u -(f_- )_v (f_+ )_u , \\ 
B & :=(f_+ )_u (f_- )_u -(f_+ )_v (f_- )_v , \\ 
C & :=(f_+ )^2_u + (f_+ )^2_v =(f_- )^2_u + (f_- )^2_v . 
\end{split} 
\label{A'} 
\end{equation} 
Then $A$, $B$, $C$ satisfy $A^2 +B^2 =C^2$. 
By \eqref{fpm}, \eqref{X+Y-X-Y+} is equivalent to 
\begin{equation} 
AA' =(f_+ )^2_v (f_- )^2_u -(f_- )^2_v (f_+ )^2_u \not= 0. 
\label{f+f-f+f-} 
\end{equation} 
Therefore $C$ is positive. 
We will prove 

\begin{lem}\label{lem:k+k-} 
The functions $k_+$, $k_-$ satisfy 
\begin{equation} 
Ak_- +B<0, \quad 
k_+ =\dfrac{Bk_- -A}{Ak_- +B} . 
\label{k+} 
\end{equation} 
\end{lem} 

\vspace{3mm} 

\par\noindent 
\textit{Proof} \ 
Applying \eqref{kpm} to \eqref{W=XY=Z}, we obtain 
\begin{equation} 
\begin{split} 
k_+ & = \dfrac{X^2_- +Y^2_- +X_+ X_- +Y_+ Y_-}{X_+ Y_- -X_- Y_+} , \\ 
k_- & =-\dfrac{X^2_+ +Y^2_+ +X_+ X_- +Y_+ Y_-}{X_+ Y_- -X_- Y_+} . 
\end{split} 
\label{k+k-} 
\end{equation} 
Applying \eqref{fpm} to \eqref{k+k-}, we obtain 
\begin{equation} 
Ak_+ -B= C\sqrt{\dfrac{1+k^2_+}{1+k^2_-}} , \quad 
Ak_- +B=-C\sqrt{\dfrac{1+k^2_-}{1+k^2_+}} , 
\label{AkpmC} 
\end{equation} 
From \eqref{AkpmC}, we obtain Lemma~\ref{lem:k+k-}. 
\hfill 
$\square$ 

\vspace{3mm} 

From \eqref{k+}, we obtain 
\begin{equation} 
1+k^2_+ =C^2 \dfrac{1+k^2_-}{(Ak_- +B)^2} 
\label{1+k2} 
\end{equation} 
and 
\begin{equation} 
(k_+ )_a 
=(AB_a -BA_a )\dfrac{1+k^2_-}{(Ak_- +B)^2} 
+ C^2 \dfrac{(k_- )_a}{(Ak_- +B)^2} 
\label{k+a} 
\end{equation} 
for $a=u, v$. 

From the first two equations in \eqref{codazzikXY}, 
we obtain 
\begin{equation} 
\begin{split} 
&  \left[ 
   \begin{array}{cc} 
    k_+ X_+ & -k_+ Y_+ \\ 
    k_- X_- &  k_- Y_- 
     \end{array} 
   \right] 
   \left[ 
   \begin{array}{c} 
    \gamma_u \\ 
    \gamma_v 
     \end{array} 
   \right] \\ 
& =\left[ 
   \begin{array}{c} 
    (X_+ )_u -(Y_+ )_v \\ 
    (X_- )_u +(Y_- )_v 
     \end{array} 
   \right]  
  +\left[ 
   \begin{array}{cc} 
   -k_+ Y_+ & -k_+ X_+ \\ 
   -k_- Y_- &  k_- X_- 
     \end{array} 
   \right] 
   \left[ 
   \begin{array}{c} 
    \lambda_u \\ 
    \lambda_v 
     \end{array} 
   \right] . 
\end{split} 
\label{gamma} 
\end{equation} 
Then noticing 
\begin{equation*} 
 k_{\pm} \not= 0, \quad 
\dfrac{\partial (f_+ , f_- )}{\partial (u, v)} =-A' \not= 0 
\end{equation*} 
and applying \eqref{fpm}, \eqref{1+k2} and \eqref{k+a} 
to \eqref{gamma}, we obtain 
\begin{equation} 
\begin{split} 
   \left[ 
   \begin{array}{c} 
    \gamma_u \\ 
    \gamma_v 
     \end{array} 
   \right] =
& -\dfrac{1}{1+k^2_-} 
   \left[ 
   \begin{array}{c} 
    (k_- )_u \\ 
    (k_- )_v 
     \end{array} 
   \right] 
  -\dfrac{A^2}{C^2} 
   \dfrac{\dfrac{\partial (f_+ , B/A )}{\partial (u, v)}}{
          \dfrac{\partial (f_+ , f_- )}{\partial (u, v)}} 
   \left[ 
   \begin{array}{c} 
    (f_- )_u \\ 
    (f_- )_v 
     \end{array} 
   \right] \\ 
& -\dfrac{1}{A'} 
   \left[ 
   \begin{array}{cc} 
    2(f_+ )_u (f_- )_u &   A \\ 
      A                & 2(f_+ )_v (f_- )_v 
     \end{array} 
   \right] 
   \left[ 
   \begin{array}{c} 
    \lambda_u \\ 
    \lambda_v 
     \end{array} 
   \right] . 
\end{split} 
\label{dgamma} 
\end{equation} 
Noticing \eqref{grad}, 
we define a function $\psi$ by 
\begin{equation} 
 \left[ \begin{array}{c} 
         (f_- )_u \\ 
         (f_- )_v 
          \end{array} 
 \right] 
=\left[ \begin{array}{cc} 
         \cos \psi & -\sin \psi \\ 
         \sin \psi &  \cos \psi 
          \end{array} 
 \right] 
 \left[ \begin{array}{c} 
         (f_+ )_v \\ 
         (f_+ )_u 
          \end{array} 
 \right] . 
\label{psi} 
\end{equation} 
We have $A=C\cos \psi$ and $B=-C\sin \psi$. 
Therefore $B/A =-\tan \psi$. 
Then \eqref{dgamma} is represented as 
\begin{equation} 
\begin{split} 
   \left[ 
   \begin{array}{c} 
    \gamma_u \\ 
    \gamma_v 
     \end{array} 
   \right] = 
& -\left[ 
 \begin{array}{c} 
  (\theta_- )_u \\ 
  (\theta_- )_v 
   \end{array} 
 \right] 
+\dfrac{\dfrac{\partial (f_+ , \psi )}{\partial (u, v)}}{
        \dfrac{\partial (f_+ , f_-  )}{\partial (u, v)}}  
 \left[ 
 \begin{array}{c} 
  (f_- )_u \\ 
  (f_- )_v 
   \end{array} 
 \right] \\ 
& -\dfrac{1}{A'} 
   \left[ 
   \begin{array}{cc} 
    2(f_+ )_u (f_- )_u &   A \\ 
      A                & 2(f_+ )_v (f_- )_v 
     \end{array} 
   \right] 
   \left[ 
   \begin{array}{c} 
    \lambda_u \\ 
    \lambda_v 
     \end{array} 
   \right] , 
\end{split} 
\label{dgamma1_2} 
\end{equation} 
where $\theta_-$ is a function given by $\tan \theta_- =k_-$. 

Let $\lambda$ be a function 
satisfying $\lambda_{uu} +\lambda_{vv} +L_0 e^{2\lambda} =0$. 
Let $\gamma$, $\theta_-$, $\psi$, $f_+$, $f_-$ be functions 
satisfying $\theta_- \not= 0$, \eqref{f+f-f+f-}, \eqref{psi} 
and \eqref{dgamma1_2}. 
We set $k_- :=\tan \theta_-$. 
Then the above functions satisfy \eqref{dgamma}. 
Let $k_+$ be as in \eqref{k+}. 
We define $W_{\pm}$, $X_{\pm}$, $Y_{\pm}$, $Z_{\pm}$ by 
\eqref{kpm} and \eqref{fpm}. 
Then these functions satisfy \eqref{W=XY=Z} and \eqref{gaussricci} 
with $\mu_1 =\gamma_u$ and $\mu_2 =\gamma_v$. 
In addition, they satisfy \eqref{codazziWXYZ} and \eqref{X+Y-X-Y+}. 
Noticing \eqref{W=XY=Z}, 
we define functions $\alpha_k$, $\beta_k$ ($k=1, 2, 3$) by 
\begin{equation} 
\alpha_1 :=\dfrac{1}{2} (Y_+ -Y_- ), \quad 
\alpha_2 :=\dfrac{1}{2} (X_+ +X_- ), \quad 
\alpha_3 :=\dfrac{1}{2} (Z_+ -Z_- ) 
\label{alpha123} 
\end{equation} 
and 
\begin{equation} 
\beta_1 :=\dfrac{1}{2} (W_+ -W_- ), \quad 
\beta_2 :=\dfrac{1}{2} (Y_+ +Y_- ), \quad 
\beta_3 :=\dfrac{1}{2} (X_+ -X_- ). 
\label{beta123} 
\end{equation} 
Then $\alpha_k$, $\beta_k$ ($k=1, 2, 3$) satisfy $S_v -T_u =ST-TS$ 
for $S$, $T$ in \eqref{ST} with $\mu_1 =\gamma_u$ and $\mu_2 =\gamma_v$. 
Therefore there exists an immersion $F$ of a neighborhood of a point 
of the $uv$-plane into $N$ with $K\equiv L_0$, $R^{\perp} \equiv 0$ 
such that the second fundamental form of $F$ does not satisfy 
the linearly dependent condition. 

Hence we obtain 

\begin{thm}\label{thm:notldL0} 
Let $F:M\longrightarrow N$ be a conformal immersion 
satisfying $K\equiv L_0$ and $R^{\perp} \equiv 0$. 
Then $F$ satisfies $k_{\pm} \not= 0$ and \eqref{X+Y-X-Y+} if and only if 
$\lambda$ with $\lambda_{uu} +\lambda_{vv} +L_0 e^{2\lambda} =0$, 
$\gamma$, $\theta_-$, $\psi$, $f_+$, $f_-$ satisfy 
$\theta_- \not= 0$, \eqref{f+f-f+f-}, \eqref{psi} and \eqref{dgamma1_2}. 
In addition, 
for functions $\lambda$, $\gamma$, $\theta_-$, $\psi$, $f_+$, $f_-$ 
as above defined on a simply connected open set $O$ 
of the $uv$-plane, 
there exists a conformal immersion $F$ of $O$ into $N$ 
with $K\equiv L_0$, $R^{\perp} \equiv 0$ 
such that the second fundamental form does not satisfy 
the linearly dependent condition, 
which is unique up to an isometry of $N$. 
\end{thm} 

In the following, suppose $L_0 =0$. 
Then we can suppose $\lambda \equiv 1$. 
Then \eqref{dgamma1_2} is represented as 
\begin{equation} 
 \left[ 
 \begin{array}{c} 
  \gamma_u \\ 
  \gamma_v 
   \end{array} 
 \right] =
-\left[ 
 \begin{array}{c} 
  (\theta_- )_u \\ 
  (\theta_- )_v 
   \end{array} 
 \right] 
+\dfrac{\dfrac{\partial (f_+ , \psi )}{\partial (u, v)}}{
        \dfrac{\partial (f_+ , f_-  )}{\partial (u, v)}}  
 \left[ 
 \begin{array}{c} 
  (f_- )_u \\ 
  (f_- )_v 
   \end{array} 
 \right] . 
\label{dgamma2} 
\end{equation} 
Since $\gamma_{uv} =\gamma_{vu}$, we obtain 
\begin{equation} 
d\left( \dfrac{\dfrac{\partial (f_+ , \psi )}{\partial (u, v)}}{
               \dfrac{\partial (f_+ , f_-  )}{\partial (u, v)}} 
 \right) 
 \wedge df_-  
=0. 
\label{gammauv} 
\end{equation} 
Therefore there exists a function $\Tilde{\xi}$ of one variable 
satisfying 
\begin{equation}
 \dfrac{\partial (f_+ , \psi )}{\partial (u, v)}
=\Tilde{\xi} (f_- )
 \dfrac{\partial (f_+ , f_-  )}{\partial (u, v)} . 
\label{gammauv2}
\end{equation}
Let $\Tilde{\Xi}$ be a primitive function of $\Tilde{\xi}$, 
i.e., a function of one variable satisfying $\Tilde{\Xi}'=\Tilde{\xi}$. 
Then by \eqref{dgamma2}, $\gamma +\theta_- -\Tilde{\Xi} (f_- )$ is constant. 
From \eqref{psi}, $(f_+ )_{uv} =(f_+ )_{vu}$ and \eqref{gammauv2}, 
we obtain 
\begin{equation} 
 \left[ 
 \begin{array}{c} 
  \psi_u \\ 
  \psi_v 
   \end{array} 
 \right] 
=(\cos^2 \psi         )\mbox{\boldmath{$a$}} 
+(\cos \psi \sin \psi )\mbox{\boldmath{$b$}}
+(\sin^2 \psi         )\mbox{\boldmath{$c$}} , 
\label{gammauv3}
\end{equation} 
where 
\begin{equation} 
\begin{split} 
  \mbox{\boldmath{$a$}} := 
& \dfrac{1}{C} 
    \left[ 
    \begin{array}{cc} 
     (f_- )_u & (f_- )_v \\ 
    -(f_- )_v & (f_- )_u 
      \end{array}
    \right] 
    \left[ 
    \begin{array}{c} 
     \Tilde{\xi} (f_- )((f_- )^2_u  -(f_- )^2_v ) \\ 
                       -(f_- )_{uu} +(f_- )_{vv} 
      \end{array}
    \right] , \\ 
  \mbox{\boldmath{$b$}} :=
& \dfrac{1}{C} 
  \left( 
   2\left[ 
    \begin{array}{cc} 
     (f_- )_u & (f_- )_v \\ 
    -(f_- )_v & (f_- )_u 
      \end{array}
    \right] 
    \left[ 
    \begin{array}{c} 
     \Tilde{\xi} (f_- )(f_- )_u (f_- )_v \\ 
                -(f_- )_{uv}  
      \end{array}
    \right] 
  \right. \\ 
& \quad \quad 
  \left.   
   +\left[ 
    \begin{array}{cc} 
     (f_- )_v & -(f_- )_u \\ 
     (f_- )_u &  (f_- )_v 
      \end{array}
    \right] 
    \left[ 
    \begin{array}{c} 
     \Tilde{\xi} (f_- )((f_- )^2_u  -(f_- )^2_v ) \\ 
                       -(f_- )_{uu} +(f_- )_{vv} 
      \end{array}
    \right] 
    \right) , \\ 
  \mbox{\boldmath{$c$}} := 
& \dfrac{2}{C} 
    \left[ 
    \begin{array}{cc} 
     (f_- )_v & -(f_- )_u \\ 
     (f_- )_u &  (f_- )_v 
      \end{array}
    \right] 
    \left[ 
    \begin{array}{c} 
     \Tilde{\xi} (f_- )(f_- )_u (f_- )_v \\ 
                -(f_- )_{uv}  
      \end{array}
    \right] . 
\end{split}
\label{abc}
\end{equation}
We can represent \eqref{gammauv3} as 
an over-determined system of polynomial type with degree two 
(see \cite{ando7} for over-determined systems of polynomial type): 
\begin{equation} 
dt=\omega_0 +t\omega_1 +t^2 \omega_2 , 
\label{abc2}
\end{equation}
where $t:=\tan \psi$ 
and $\omega_l$ ($l=0, 1, 2$) are $1$-forms given by 
\begin{equation}
\omega_0 =[du \ dv]\mbox{\boldmath{$a$}} , \quad 
\omega_1 =[du \ dv]\mbox{\boldmath{$b$}} , \quad 
\omega_2 =[du \ dv]\mbox{\boldmath{$c$}} . 
\label{omega} 
\end{equation}
By $d(dt)=0$, we obtain 
\begin{equation*} 
\begin{split}
&  \Omega_0 +t\Omega_1 +t^2 \Omega_2 =0 \\ 
& (\Omega_0 :=d\omega_0 + \omega_0 \wedge \omega_1  , \quad 
   \Omega_1 :=d\omega_1 +2\omega_0 \wedge \omega_2  , \quad 
   \Omega_2 :=d\omega_2 +2\omega_1 \wedge \omega_2 ).  
\end{split} 
\label{Omega}   
\end{equation*} 

Let $\Tilde{\xi}$ be a function of one variable. 
Let $f_-$ be a function of two variables $u$, $v$ 
satisfying $(f_- )^2_u +(f_- )^2_v \not= 0$. 
If $f_-$ satisfies $(\Omega_0 , \Omega_1 , \Omega_2 )=(0, 0, 0)$, 
then for an arbitrarily given initial value, 
\eqref{abc2} with \eqref{abc} and \eqref{omega} has a unique solution $t$; 
if $f_-$ does not satisfy $(\Omega_0 , \Omega_1 , \Omega_2 )=(0, 0, 0)$, 
then the number of solutions of \eqref{abc2} is zero, one or two 
(\cite{ando7}). 
Let $t$ be a solution of \eqref{abc2} satisfying 
\begin{equation} 
2(f_- )_u (f_- )_v t +(f_-)^2_u -(f_- )^2_v \not= 0. 
\label{t} 
\end{equation} 
Let $\psi$ be a function given by $t=\tan \psi$. 
Then there exists a function $f_+$ 
satisfying $A\not= 0$, \eqref{psi} and \eqref{gammauv2}. 
Since $f_-$ and $t$ satisfy \eqref{t}, we obtain $A'\not= 0$. 
Therefore there exist functions $\gamma$, $\theta_-$ 
of two variables $u$, $v$ satisfying $\theta_- \not= 0$ 
and \eqref{dgamma2}. 
Then $\lambda \equiv 1$, $\gamma$, $\theta_-$, $\psi$, $f_+$, $f_-$ 
satisfy \eqref{f+f-f+f-} and \eqref{dgamma1_2}. 

Therefore from Theorem~\ref{thm:notldL0}, we obtain 

\begin{thm}\label{thm:notld} 
Let $F:M\longrightarrow E^4$ be a conformal immersion 
satisfying $K\equiv 0$ and $R^{\perp} \equiv 0$. 
Then $F$ satisfies $k_{\pm} \not= 0$ and \eqref{X+Y-X-Y+} if and only if 
the functions $\gamma$, $\theta_-$, $f_-$, $t=\tan \psi$ satisfy 
$\theta_- \not= 0$, $(f_- )^2_u +(f_- )^2_v \not= 0$, 
\eqref{dgamma2}, \eqref{abc2} with \eqref{abc} and \eqref{omega} 
for a function $\Tilde{\xi}$ of one variable and \eqref{t}. 
In addition, 
for functions $f_-$, $\psi$ as above 
defined on a simply connected open set $O$ of the $uv$-plane, 
there exists a conformal immersion $F$ of $O$ into $E^4$ 
with $K\equiv 0$, $R^{\perp} \equiv 0$ 
such that the second fundamental form does not satisfy 
the linearly dependent condition. 
\end{thm} 

\begin{ex} 
Suppose that $f_+$, $\psi$ satisfy $d\psi \wedge df_+ =0$, 
which is equivalent to $\Tilde{\xi} \equiv 0$. 
Then we have \eqref{gammauv}, 
and there exists a function $\xi$ of one variable 
satisfying $\psi =\xi (f_+ )$. 
From \eqref{psi} and $(f_- )_{uv} =(f_- )_{vu}$, we obtain 
\begin{equation} 
\begin{split} 
&   c_+ ((f_+ )_{uu} -(f_+ )_{vv} ) 
  +2s_+  (f_+ )_{uv} \\ 
& +\xi'(f_+ )(- s_+ ((f_+ )^2_u -(f_+ )^2_v ) 
              +2c_+  (f_+ )_u (f_+ )_v )=0, 
\end{split} 
\label{xieq} 
\end{equation} 
where 
\begin{equation} 
c_+ :=\cos (\xi (f_+ )), \quad 
s_+ :=\sin (\xi (f_+ )). 
\label{c+s+} 
\end{equation} 
In particular, 
if $\psi$ is constant, 
then we have 
\begin{equation*} 
  (\cos \psi )((f_+ )_{uu} -(f_+ )_{vv} ) 
+2(\sin \psi ) (f_+ )_{uv} =0. 
\end{equation*} 
Therefore for a function $\xi$ of one variable 
and a function $f_+$ of two variables $u$, $v$ satisfying \eqref{xieq} 
with \eqref{c+s+}, $(f_+ )^2_u +(f_+ )^2_v >0$ 
and $\cos \xi (f_+ )\not= 0$, 
there exists a function $f_-$ satisfying \eqref{psi} 
with $\psi :=\xi (f_+ )$. 
Then we have $A\not= 0$. 
In addition, 
if we suppose $c_+ ((f_+ )^2_u -(f_+ )^2_v )+2s_+ (f_+ )_u (f_+ )_v 
              \not= 0$, 
then $A'\not= 0$. 
\end{ex} 

\section{Space-like surfaces with flat normal connection 
in neutral 4-dimensional space forms}\label{sect:sn} 

\setcounter{equation}{0} 

Let $N$ be a 4-dimensional neutral space form 
with constant sectional curvature $L_0$. 
If $L_0 =0$, 
then we can suppose $N=E^4_2$; 
if $L_0 >0$, 
then we can suppose $N=\{ x\in E^5_2 \ | \ \langle x, x\rangle =1/L_0 \}$; 
if $L_0 <0$, 
then we can suppose $N=\{ x\in E^5_3 \ | \ \langle x, x\rangle =1/L_0 \}$, 
where $\langle \ , \ \rangle$ is the metric of $E^5_2$ or $E^5_3$.  
Let $h$ be the metric of $N$. 
Let $M$ be a Riemann surface 
and $F:M\longrightarrow N$ a space-like and conformal immersion 
of $M$ into $N$. 
Let $(u, v)$ be as in Section~\ref{sect:R}. 
Then the induced metric $g$ on $M$ by $F$ is represented 
as $g=e^{2\lambda} (du^2 +dv^2 )$ 
for a real-valued function $\lambda$. 
Let $\Tilde{\nabla}$ denote the Levi-Civita connection 
of $E^4_2$, $E^5_2$ or $E^5_3$ according to $L_0 =0$, $>0$ or $<0$. 
We set $T_1 :=dF(\partial_u )$, $T_2 :=dF(\partial_v )$. 
Let $N_1$, $N_2$ be normal vector fields of $F$ satisfying 
$$h(N_1 , N_1 )=h(N_2 , N_2 )=-e^{2\lambda} , \quad 
  h(N_1 , N_2 )=0.$$ 
Then we have 
\begin{equation} 
\begin{split} 
\Tilde{\nabla}_{T_1} (T_1 \ T_2 \ N_1 \ N_2 \ F) 
& =                  (T_1 \ T_2 \ N_1 \ N_2 \ F)S, \\ 
\Tilde{\nabla}_{T_2} (T_1 \ T_2 \ N_1 \ N_2 \ F) 
& =                  (T_1 \ T_2 \ N_1 \ N_2 \ F)T, 
\end{split} 
\label{dt1t2sn} 
\end{equation} 
where 
\begin{equation*} 
\begin{split} 
S & =\left[ \begin{array}{ccccc} 
         \lambda_u         & \lambda_v &  \alpha_1  &  \beta_1   & 1 \\ 
        -\lambda_v         & \lambda_u &  \alpha_2  &  \beta_2   & 0 \\ 
         \alpha_1          & \alpha_2  &  \lambda_u & -\mu_1     & 0 \\ 
         \beta_1           & \beta_2   &  \mu_1     &  \lambda_u & 0 \\ 
         -L_0 e^{2\lambda} &  0        &   0        &   0        & 0 
              \end{array} 
     \right] , \\ 
T & =\left[ \begin{array}{ccccc} 
         \lambda_v & -\lambda_u        &  \alpha_2  &  \beta_2   & 0 \\ 
         \lambda_u &  \lambda_v        &  \alpha_3  &  \beta_3   & 1 \\ 
         \alpha_2  &  \alpha_3         &  \lambda_v & -\mu_2     & 0 \\ 
         \beta_2   &  \beta_3          &  \mu_2     &  \lambda_v & 0 \\ 
          0        & -L_0 e^{2\lambda} &   0        &   0        & 0 
              \end{array} 
    \right] , 
\end{split} 
\label{STsn} 
\end{equation*} 
and $\alpha_k$, $\beta_k$ ($k=1, 2, 3$) and $\mu_l$ ($l=1, 2$) are 
real-valued functions. 
From \eqref{dt1t2sn}, we obtain $S_v -T_u =ST-TS$. 
This is equivalent to 
the system of the equations of Gauss, Codazzi and Ricci. 
The equation of Gauss is given by 
\begin{equation*} 
  \lambda_{uu} +\lambda_{vv} +L_0 e^{2\lambda} 
= \alpha_1 \alpha_3 +\beta_1 \beta_3 -\alpha^2_2 -\beta^2_2 . 
\label{gausssn}
\end{equation*}
The equations of Codazzi are given by \eqref{codazzi}. 
The equation of Ricci is given by 
\begin{equation*} 
  (\mu_1 )_v -(\mu_2 )_u =-\alpha_1 \beta_2 +\alpha_2 \beta_1 
                          -\alpha_2 \beta_3 +\alpha_3 \beta_2 . 
\label{riccisn} 
\end{equation*}

Theorem~\ref{thm:pnvf} holds 
for a space-like and conformal immersion $F:M\longrightarrow N$, 
and the proof is exactly given by the proof of Theorem~\ref{thm:pnvf} 
in Subsection~\ref{subsect:Rpnvf}.  

We have an analogue of Theorem~\ref{thm:ldnotpnvf} 
for space-like surfaces in $N$. 
In addition, 
we have analogues of Theorems~\ref{thm:notldL0} and \ref{thm:notld} 
for space-like surfaces in $N$. 

\section{Time-like surfaces with flat normal connection 
in neutral 4-dimensional space forms}\label{sect:tn} 

\setcounter{equation}{0} 

\subsection{The Gauss-Codazzi-Ricci equations}\label{subsect:tnGCR} 

Let $N$, $h$ be as in Section~\ref{sect:sn}. 
Let $M$ be a Lorentz surface 
and $F:M\longrightarrow N$ a time-like and conformal immersion 
of $M$ into $N$. 
Let $(u, v)$ be local coordinates of $M$ 
which are compatible with the paracomplex structure of $M$. 
Then the induced metric $g$ on $M$ by $F$ is represented 
as $g=e^{2\lambda} (du^2 -dv^2 )$ 
for a real-valued function $\lambda$. 
Let $\Tilde{\nabla}$, $T_1$, $T_2$ be as in Section~\ref{sect:sn}. 
Let $N_1$, $N_2$ be normal vector fields of $F$ satisfying 
$$h(N_1 , N_1 )=-h(N_2 , N_2 )=e^{2\lambda} , \quad 
  h(N_1 , N_2 )=0.$$ 
Then we have 
\begin{equation} 
\begin{split} 
\Tilde{\nabla}_{T_1} (T_1 \ T_2 \ N_1 \ N_2 \ F) 
& =                  (T_1 \ T_2 \ N_1 \ N_2 \ F)S, \\ 
\Tilde{\nabla}_{T_2} (T_1 \ T_2 \ N_1 \ N_2 \ F) 
& =                  (T_1 \ T_2 \ N_1 \ N_2 \ F)T, 
\end{split} 
\label{dt1t2tn} 
\end{equation} 
where 
\begin{equation} 
\begin{split} 
S & =\left[ \begin{array}{ccccc} 
         \lambda_u         & \lambda_v & -\alpha_1  &  \beta_1   & 1 \\ 
         \lambda_v         & \lambda_u &  \alpha_2  & -\beta_2   & 0 \\ 
         \alpha_1          & \alpha_2  &  \lambda_u &  \mu_1     & 0 \\ 
         \beta_1           & \beta_2   &  \mu_1     &  \lambda_u & 0 \\ 
         -L_0 e^{2\lambda} &  0        &   0        &   0        & 0 
              \end{array} 
     \right] , \\ 
T & =\left[ \begin{array}{ccccc} 
         \lambda_v &  \lambda_u        & -\alpha_2  &  \beta_2   & 0 \\ 
         \lambda_u &  \lambda_v        &  \alpha_3  & -\beta_3   & 1 \\ 
         \alpha_2  &  \alpha_3         &  \lambda_v &  \mu_2     & 0 \\ 
         \beta_2   &  \beta_3          &  \mu_2     &  \lambda_v & 0 \\ 
          0        &  L_0 e^{2\lambda} &   0        &   0        & 0 
              \end{array} 
    \right] , 
\end{split} 
\label{STtn} 
\end{equation} 
and $\alpha_k$, $\beta_k$ ($k=1, 2, 3$) and $\mu_l$ ($l=1, 2$) are 
real-valued functions. 
From \eqref{dt1t2tn}, we obtain $S_v -T_u =ST-TS$. 
This is equivalent to 
the system of the equations of Gauss, Codazzi and Ricci. 
The equation of Gauss is given by 
\begin{equation} 
  \lambda_{uu} -\lambda_{vv} +L_0 e^{2\lambda} 
= \alpha_1 \alpha_3 -\beta_1 \beta_3 -\alpha^2_2 +\beta^2_2 . 
\label{gausstn}
\end{equation}
The equations of Codazzi are given by 
\begin{equation} 
\begin{split}
  (\alpha_1 )_v -(\alpha_2 )_u 
& =\alpha_2 \lambda_u -\alpha_3 \lambda_v 
   +\beta_2 \mu_1     -\beta_1  \mu_2  , \\ 
   (\alpha_2 )_v -(\alpha_3 )_u 
& = \alpha_1 \lambda_u -\alpha_2 \lambda_v 
   +\beta_3  \mu_1     -\beta_2  \mu_2  , \\ 
   (\beta_1 )_v -(\beta_2 )_u 
& = \beta_2  \lambda_u -\beta_3  \lambda_v 
   +\alpha_2 \mu_1     -\alpha_1 \mu_2  , \\ 
   (\beta_2 )_v -(\beta_3 )_u 
& = \beta_1  \lambda_u -\beta_2  \lambda_v 
   +\alpha_3 \mu_1     -\alpha_2 \mu_2 .  
\end{split} 
\label{codazzitn} 
\end{equation} 
The equation of Ricci is given by 
\begin{equation} 
  (\mu_1 )_v -(\mu_2 )_u = \alpha_1 \beta_2 -\alpha_2 \beta_1 
                          -\alpha_2 \beta_3 +\alpha_3 \beta_2 . 
\label{riccitn} 
\end{equation}

\subsection{Time-like surfaces 
with \mbox{\boldmath{$R^{\perp} \equiv 0$}} and 
parallel normal vector fields}\label{subsect:tnpnvf} 

Suppose that $\nabla^{\perp}$ is flat. 
Then there exists a function $\gamma$ satisfying 
$\gamma_u =\mu_1$, $\gamma_v =\mu_2$. 
In addition, we obtain 
\begin{equation*} 
 \alpha_1 \beta_2 -\alpha_2 \beta_1 
-\alpha_2 \beta_3 +\alpha_3 \beta_2 =0. 
\label{riccifnctn} 
\end{equation*} 
We say that the second fundamental form $\sigma$ of $F$ satisfies 
\begin{itemize} 
\item{the \textit{space-like linearly dependent condition\/} 
if $F$ satisfies 
\begin{equation*} 
 {\rm cosh}\,t_+ (\alpha_1 , \alpha_2 , \alpha_3 ) 
+{\rm sinh}\,t_+ (\beta_1  , \beta_2  , \beta_3  ) 
=(0, 0, 0) 
\label{tsp} 
\end{equation*} 
for a function $t_+$ of $u$, $v$;} 
\item{the \textit{time-like linearly dependent condition\/} 
if $F$ satisfies 
\begin{equation*} 
 {\rm sinh}\,t_- (\alpha_1 , \alpha_2 , \alpha_3 ) 
+{\rm cosh}\,t_- (\beta_1  , \beta_2  , \beta_3  ) 
=(0, 0, 0) 
\label{ttm} 
\end{equation*} 
for a function $t_-$ of $u$, $v$;} 
\item{the \textit{light-like linearly dependent condition\/} 
if $F$ satisfies 
\begin{equation*} 
             (\alpha_1 , \alpha_2 , \alpha_3 ) 
+\varepsilon (\beta_1  , \beta_2  , \beta_3  ) 
=(0,0, 0) 
\label{el} 
\end{equation*} 
for $\varepsilon \in \{ 1, -1\}$.} 
\end{itemize} 
We say that the second fundamental form $\sigma$ of $F$ satisfies 
the \textit{linearly dependent condition\/} 
if $\sigma$ satisfies 
the space-like, time-like or light-like linearly dependent condition. 
If $F$ has a light-like parallel normal vector field, 
then $\sigma$ satisfies the light-like linearly dependent condition. 
If $\sigma$ satisfies the light-like linearly dependent condition, 
then by \eqref{gausstn}, we have $K\equiv L_0$. 

Referring to the proof of Theorem~\ref{thm:pnvf}, 
we can prove 

\begin{thm}\label{thm:pnvftn} 
Suppose that $\nabla^{\perp}$ is flat. 
\begin{itemize} 
\item[{\rm (a)}]{Suppose $K\not= L_0$. 
Then on a neighborhood of each point of $M$, 
$F$ has a space-like $($respectively, time-like$)$ parallel normal 
vector field 
if and only if the second fundamental form $\sigma$ of $F$ satisfies 
the space-like $($respectively, time-like$)$ linearly dependent condition. 
In addition, if one of these conditions holds, 
then $\gamma -t_+ ($respectively, $\gamma -t_- )$ is constant.} 
\item[{\rm (b)}]{Suppose $K\equiv L_0$. 
Then on a neighborhood of each point of $M$, 
$F$ has 
a space-like $($respectively, time-like$)$ parallel normal vector field 
if and only if $\sigma$ satisfies 
the space-like $($respectively, time-like$)$ linearly dependent condition 
so that $\gamma -t_+ ($respectively, $\gamma -t_- )$ is constant.} 
\end{itemize} 
If $F$ has 
a space-like $($respectively, time-like$)$ parallel normal vector field, 
then it is given by 
$e^{-\lambda} ({\rm cosh}\,t_+ N_1 -{\rm sinh}\,t_+ N_2 )$ 
$($respectively, 
$e^{-\lambda} ({\rm sinh}\,t_- N_1 -{\rm cosh}\,t_- N_2 ))$.  
\end{thm}

\begin{thm}\label{thm:pnvftnll} 
Suppose that $\nabla^{\perp}$ is flat and $K\equiv L_0$. 
Then on a neighborhood of each point of $M$, 
$F$ has a light-like parallel normal vector field 
if and only if the second fundamental form $\sigma$ of $F$ satisfies 
the light-like linearly dependent condition. 
In addition, if one of these conditions holds, 
then a light-like parallel normal vector field is given by 
$ae^{-\lambda +\varepsilon \gamma} (N_1 -\varepsilon N_2 )$ 
for a nonzero constant $a$ and $\varepsilon \in \{ 1, -1\}$. 
\end{thm} 

\begin{rem}\label{rem:K=L_0} 
Suppose that $F$ has zero mean curvature vector. 
Then the following are mutually equivalent (\cite{ando5}, \cite{ando9}): 
\begin{itemize} 
\item[{\rm (a)}]{$K\equiv L_0$;} 
\item[{\rm (b)}]{$\nabla^{\perp}$ is flat, and 
the paraholomorphic quartic differential $Q$ is zero or null;} 
\item[{\rm (c)}]{the covariant derivatives of 
both of the time-like twistor lifts of $F$ are zero or light-like;} 
\item[{\rm (d)}]{$F$ satisfies one of the following: 
\begin{itemize} 
\item[{\rm (d-i)}]{a light-like normal vector field of $F$ is parallel;} 
\item[{\rm (d-ii)}]{the shape operator of any normal vector field is 
zero or light-like.} 
\end{itemize}} 
\end{itemize} 
We can refer to \cite{ando4}, \cite{ando5}, \cite{ando9}, \cite{ando8}, 
\cite{HM}, \cite{JR} for the time-like twistor lifts. 
The conformal Gauss maps of time-like surfaces in 
$3$-dimensional Lorentzian space forms of Willmore type 
with vanishing paraholomorphic quartic differential satisfy (d-i) 
and (a) with $L_0 =1$ (\cite{ando4}). 
In particular, 
the conformal Gauss maps of time-like surfaces in 
the Lorentz-Minkowski $3$-space $E^3_1$ with zero mean curvature 
satisfy (d-i) (we can refer to \cite{AHHK} for such surfaces).  
If $F$ satisfies (d-ii), 
then the second fundamental form $\sigma$ of $F$ satisfies 
the linearly dependent condition. 
See \cite{ando5} for a characterization of time-like surfaces in $N$ 
with zero mean curvature vector and (d-ii). 
There exist time-like surfaces in $N$ 
with zero mean curvature vector, 
$K\not= L_0$ and $R^{\perp} \equiv 0$ (\cite{ando9}), 
and by Theorem~\ref{thm:pnvftn}, 
such surfaces have space-like or time-like parallel normal vector fields. 
\end{rem} 

\subsection{Time-like surfaces with \mbox{\boldmath{$R^{\perp} \equiv 0$}} 
and \mbox{\boldmath{$K\equiv L_0$}} 
which do not admit any parallel normal vector fields}\label{subsect:tnnpnvf} 

We have two analogues of Theorem~\ref{thm:ldnotpnvf} 
for time-like surfaces in $N$, 
according to the space-like linearly dependent condition  
and          the time-like  linearly dependent condition. 

In the following, 
we study a time-like and conformal immersion $F:M\longrightarrow N$ 
of $M$ into $N$ with $R^{\perp} \equiv 0$ and $K\equiv L_0$ 
such that the second fundamental form of $F$ does not satisfy 
the linearly dependent condition. 

Let $W_{\pm}$, $X_{\pm}$, $Y_{\pm}$, $Z_{\pm}$ be as in \eqref{WXYZ}. 
Then we have \eqref{W=XY=Z}.  
By \eqref{gausstn} and \eqref{riccitn}, we obtain 
\begin{equation} 
-W_{\pm} X_{\pm} +Y_{\pm} Z_{\pm} 
=\lambda_{uu} -\lambda_{vv} +L_0 e^{2\lambda} 
 \pm ((\mu_1 )_v -(\mu_2 )_u ). 
\label{gaussriccitn} 
\end{equation} 
We can rewrite \eqref{codazzitn} into 
\begin{equation} 
\begin{split}
\pm  (Y_{\pm} )_v -(X_{\pm} )_u 
  & = W_{\pm} (\lambda_u \mp \mu_2 )+Z_{\pm} (\mp \lambda_v +\mu_1 ), \\ 
\pm  (W_{\pm} )_v -(Z_{\pm} )_u 
  & = Y_{\pm} (\lambda_u \mp \mu_2 )+X_{\pm} (\mp \lambda_v +\mu_1 ).  
\end{split} 
\label{codazziWXYZtn} 
\end{equation} 
  
Suppose $R^{\perp} \equiv 0$ and $K\equiv L_0$. 
Then from \eqref{gaussriccitn}, 
we obtain 
\begin{equation} 
W_{\pm} X_{\pm} -Y_{\pm} Z_{\pm} =0. 
\label{WXYZ0tn} 
\end{equation} 
We suppose that there exist functions $k_{\pm}$ 
satisfying $k_{\pm} \not= 0$ and 
\begin{equation} 
(W_{\pm} , Z_{\pm} )=k_{\pm} (Y_{\pm} , X_{\pm} ). 
\label{kpmtn} 
\end{equation} 
In addition, we suppose $|k_{\pm} |\not= 1$. 
We will prove 

\begin{lem}\label{lem:fXYtn} 
There exist functions $f_+$, $f_-$ satisfying 
\begin{equation} 
\begin{array}{lcl} 
 X_+ = \dfrac{(f_+ )_v}{\sqrt{|k^2_+ -1|}} , & \ &  
 Y_+ = \dfrac{(f_+ )_u}{\sqrt{|k^2_+ -1|}} ,   \\ 
 X_- =-\dfrac{(f_- )_v}{\sqrt{|k^2_- -1|}} , & \ &  
 Y_- = \dfrac{(f_- )_u}{\sqrt{|k^2_- -1|}} . 
  \end{array} 
\label{fpmtn}
\end{equation} 
\end{lem} 

\vspace{3mm} 

\par\noindent 
\textit{Proof} \ 
Applying \eqref{kpmtn} to \eqref{codazziWXYZtn}, we obtain 
\begin{equation} 
\begin{split} 
      \pm (Y_{\pm} )_v -(X_{\pm} )_u 
& =k_{\pm} Y_{\pm} (    \lambda_u \mp \gamma_v ) 
  +k_{\pm} X_{\pm} (\mp \lambda_v  +  \gamma_u ), \\  
  (k_{\pm} Y_{\pm} )_v \mp (k_{\pm} X_{\pm} )_u 
& = Y_{\pm} (\pm \lambda_u -\gamma_v )+X_{\pm} (-\lambda_v \pm \gamma_u ).  
\end{split} 
\label{codazzikXYtn} 
\end{equation} 
From \eqref{codazzikXYtn}, we obtain 
\begin{equation*} 
     \left( \sqrt{|k^2_{\pm} -1|} X_{\pm} \right)_u 
=\pm \left( \sqrt{|k^2_{\pm} -1|} Y_{\pm} \right)_v . 
\label{kXYpmtn} 
\end{equation*} 
Therefore there exist functions $f_+$, $f_-$ 
satisfying \eqref{fpmtn}. 
Hence we obtain Lemma~\ref{lem:fXYtn}. 
\hfill 
$\square$ 

\vspace{3mm} 

We will prove 

\begin{lem}\label{lem:f+2f-2tn} 
The functions $f_+$, $f_-$ satisfy 
\begin{equation} 
(f_+ )^2_u - (f_+ )^2_v =\varepsilon ((f_- )^2_u - (f_- )^2_v )  
\label{gradtn} 
\end{equation} 
for $\varepsilon \in \{ 1, -1\}$. 
\end{lem} 

\vspace{3mm} 

\par\noindent 
\textit{Proof} \ 
Since $R^{\perp} \equiv 0$ is equivalent to 
\begin{equation}  
 W^2_+ +X^2_+ +Y^2_- +Z^2_-  
=W^2_- +X^2_- +Y^2_+ +Z^2_+ , 
\label{X2Y2Z2W2tn} 
\end{equation} 
from \eqref{kpmtn}, \eqref{fpmtn} and \eqref{X2Y2Z2W2tn}, 
we obtain \eqref{gradtn}. 
\hfill 
$\square$ 

\vspace{3mm} 

Noticing $|k_{\pm} |\not= 1$ and 
referring to the proof of Lemma~\ref{lem:XYXY},  
We can prove 

\begin{lem}\label{lem:XYXYtn} 
Suppose \eqref{X+Y-X-Y+}. 
Then the second fundamental form of $F$ does not satisfy 
the linearly dependent condition. 
\end{lem} 

In the following, we suppose \eqref{X+Y-X-Y+}. 
Let $A$, $A'$ be as in \eqref{A'} and set 
\begin{equation*} 
\begin{split} 
B & := (f_+ )_u (f_- )_u +(f_+ )_v (f_- )_v , \\ 
C & := (f_+ )^2_u - (f_+ )^2_v =\varepsilon ((f_- )^2_u - (f_- )^2_v ). 
\end{split} 
\label{ABCtn} 
\end{equation*} 
Then $A$, $B$, $C$ satisfy $A^2 -B^2 =-\varepsilon C^2$. 
By \eqref{f+f-f+f-}, $C$ is nonzero. 
We will prove 

\begin{lem}\label{lem:k+k-tn} 
The functions $k_+$, $k_-$ satisfy 
\begin{equation} 
(Ak_- +B)C\not= 0, \quad 
k_+ =\dfrac{Bk_- +A}{Ak_- +B} . 
\label{k+tn} 
\end{equation} 
\end{lem} 

\vspace{3mm} 

\par\noindent 
\textit{Proof} \ 
Applying \eqref{kpmtn} to \eqref{W=XY=Z}, we obtain 
\begin{equation} 
\begin{split} 
k_+ & =-\dfrac{X^2_- -Y^2_- +X_+ X_- -Y_+ Y_-}{X_+ Y_- -X_- Y_+} , \\ 
k_- & = \dfrac{X^2_+ -Y^2_+ +X_+ X_- -Y_+ Y_-}{X_+ Y_- -X_- Y_+} . 
\end{split} 
\label{k+k-tn} 
\end{equation} 
Applying \eqref{fpmtn} to \eqref{k+k-tn}, we obtain 
\begin{equation} 
Ak_+ -B= \varepsilon C\sqrt{\dfrac{|k^2_+ -1|}{|k^2_- -1|}} , \quad 
Ak_- +B=-            C\sqrt{\dfrac{|k^2_- -1|}{|k^2_+ -1|}} . 
\label{AkpmCtn} 
\end{equation} 
From \eqref{AkpmCtn}, we obtain Lemma~\ref{lem:k+k-tn}. 
\hfill 
$\square$ 

\vspace{3mm} 

From \eqref{k+tn}, we obtain 
\begin{equation} 
k^2_+ -1=\varepsilon C^2 \dfrac{k^2_- -1}{(Ak_- +B)^2} 
\label{1+k2tn} 
\end{equation} 
and 
\begin{equation} 
(k_+ )_a 
=(AB_a -BA_a )\dfrac{k^2_- -1}{(Ak_- +B)^2} 
+\varepsilon  C^2 \dfrac{(k_- )_a}{(Ak_- +B)^2} 
\label{k+atn} 
\end{equation} 
for $a=u, v$. 

From the first two equations in \eqref{codazzikXYtn}, 
we obtain 
\begin{equation} 
\begin{split} 
&  \left[ 
   \begin{array}{cc} 
    k_+ X_+ & -k_+ Y_+ \\ 
    k_- X_- &  k_- Y_- 
     \end{array} 
   \right] 
   \left[ 
   \begin{array}{c} 
    \gamma_u \\ 
    \gamma_v 
     \end{array} 
   \right] \\ 
& =\left[ 
   \begin{array}{c} 
   -(X_+ )_u +(Y_+ )_v \\ 
   -(X_- )_u -(Y_- )_v 
     \end{array} 
   \right]  
  +\left[ 
   \begin{array}{cc} 
   -k_+ Y_+ &  k_+ X_+ \\ 
   -k_- Y_- & -k_- X_- 
     \end{array} 
   \right] 
   \left[ 
   \begin{array}{c} 
    \lambda_u \\ 
    \lambda_v 
     \end{array} 
   \right] . 
\end{split} 
\label{gammatn} 
\end{equation} 
Applying \eqref{fpmtn}, \eqref{1+k2tn} and \eqref{k+atn} 
to \eqref{gammatn}, we obtain 
\begin{equation} 
\begin{split} 
   \left[ 
   \begin{array}{c} 
    \gamma_u \\ 
    \gamma_v 
     \end{array} 
   \right] =
&  \dfrac{1}{k^2_- -1} 
   \left[ 
   \begin{array}{c} 
    (k_- )_u \\ 
    (k_- )_v 
     \end{array} 
   \right] 
  +\dfrac{A^2}{\varepsilon C^2} 
   \dfrac{\dfrac{\partial (f_+ , B/A )}{\partial (u, v)}}{
          \dfrac{\partial (f_+ , f_- )}{\partial (u, v)}} 
   \left[ 
   \begin{array}{c} 
    (f_- )_u \\ 
    (f_- )_v 
     \end{array} 
   \right] \\ 
& -\dfrac{1}{A'} 
   \left[ 
   \begin{array}{cc} 
    2(f_+ )_u (f_- )_u & -  A \\ 
      A                & -2(f_+ )_v (f_- )_v 
     \end{array} 
   \right] 
   \left[ 
   \begin{array}{c} 
    \lambda_u \\ 
    \lambda_v 
     \end{array} 
   \right] . 
\end{split} 
\label{dgammatn} 
\end{equation} 
\begin{itemize} 
\item{Suppose $\varepsilon =1$. 
Noticing \eqref{gradtn}, 
we define a function $\rho$ by 
\begin{equation} 
 \left[ \begin{array}{c} 
         (f_- )_u \\ 
         (f_- )_v 
          \end{array} 
 \right] 
=\left[ \begin{array}{cc} 
         \delta\,{\rm cosh}\,\rho &         {\rm sinh}\,\rho \\ 
                 {\rm sinh}\,\rho & \delta\,{\rm cosh}\,\rho 
          \end{array} 
 \right] 
 \left[ \begin{array}{c} 
         (f_+ )_u \\ 
        -(f_+ )_v 
          \end{array} 
 \right] , 
\label{e=+rho} 
\end{equation} 
where $\delta \in \{ 1, -1\}$. 
Then by $A\not= 0$, we have $\rho \not= 0$. 
In addition, we have $A=       C\,{\rm sinh}\,\rho$ 
and                  $B=\delta C\,{\rm cosh}\,\rho$. 
Therefore $B/A =\delta\,{\rm coth}\, \rho$.}  
\item{Suppose $\varepsilon =-1$. 
We define a function $\rho$ by 
\begin{equation} 
 \left[ \begin{array}{c} 
         (f_- )_u \\ 
         (f_- )_v 
          \end{array} 
 \right] 
=\left[ \begin{array}{cc} 
         \delta\,{\rm cosh}\,\rho &         {\rm sinh}\,\rho \\ 
                 {\rm sinh}\,\rho & \delta\,{\rm cosh}\,\rho 
          \end{array} 
 \right] 
 \left[ \begin{array}{c} 
         (f_+ )_v \\ 
        -(f_+ )_u 
          \end{array} 
 \right] 
\label{e=-rho} 
\end{equation} 
($\delta \in \{ 1, -1\}$). 
Then we have $A=-\delta C\,{\rm cosh}\,\rho$ 
and          $B=-       C\,{\rm sinh}\,\rho$. 
Therefore $B/A =\delta\,{\rm tanh}\, \rho$.} 
\end{itemize} 
Then \eqref{dgammatn} is represented as 
\begin{equation} 
\begin{split} 
   \left[ 
   \begin{array}{c} 
    \gamma_u \\ 
    \gamma_v 
     \end{array} 
   \right] =
&  \left[ 
   \begin{array}{c} 
    (t_- )_u \\ 
    (t_- )_v 
     \end{array} 
   \right] 
  -\delta 
   \dfrac{\dfrac{\partial (f_+ , \rho )}{\partial (u, v)}}{
          \dfrac{\partial (f_+ , f_-  )}{\partial (u, v)}}  
   \left[ 
   \begin{array}{c} 
    (f_- )_u \\ 
    (f_- )_v 
     \end{array} 
   \right] \\ 
& -\dfrac{1}{A'} 
   \left[ 
   \begin{array}{cc} 
    2(f_+ )_u (f_- )_u & -  A \\ 
      A                & -2(f_+ )_v (f_- )_v 
     \end{array} 
   \right] 
   \left[ 
   \begin{array}{c} 
    \lambda_u \\ 
    \lambda_v 
     \end{array} 
   \right] , 
\end{split} 
\label{dgamma1_2tn} 
\end{equation} 
where 
\begin{equation} 
t_- =\dfrac{1}{2} \log \left| \dfrac{k_- -1}{k_- +1} \right| . 
\label{t-}
\end{equation} 

Let $\lambda$ be a function 
satisfying $\lambda_{uu} -\lambda_{vv} +L_0 e^{2\lambda} =0$. 
Let $\gamma$, $t_-$, $\rho$, $f_+$, $f_-$ be functions 
satisfying $t_- \not= 0$, \eqref{f+f-f+f-}, \eqref{e=+rho} 
and \eqref{dgamma1_2tn}. 
Let $k_-$ be a function satisfying \eqref{t-}. 
Then we have 
\begin{equation*} 
k_- =\dfrac{1+\varepsilon' e^{2t_-}}{1-\varepsilon' e^{2t_-}} 
\end{equation*} 
for $\varepsilon' =1$ or $-1$. 
Then the above functions satisfy \eqref{dgammatn} 
with $\varepsilon =1$. 
Let $k_+$ be as in \eqref{k+tn}. 
We define $W_{\pm}$, $X_{\pm}$, $Y_{\pm}$, $Z_{\pm}$ by 
\eqref{kpmtn} and \eqref{fpmtn}. 
Then these functions satisfy \eqref{W=XY=Z}, 
\eqref{gaussriccitn} with $\mu_1 =\gamma_u$ and $\mu_2 =\gamma_v$. 
In addition, they satisfy \eqref{codazziWXYZtn} and \eqref{X+Y-X-Y+}. 
Noticing \eqref{W=XY=Z}, 
we define functions $\alpha_k$, $\beta_k$ ($k=1, 2, 3$) 
by \eqref{alpha123} and \eqref{beta123}. 
Then $\alpha_k$, $\beta_k$ ($k=1, 2, 3$) satisfy $S_v -T_u =ST-TS$ 
for $S$, $T$ in \eqref{STtn} 
with $\mu_1 =\gamma_u$ and $\mu_2 =\gamma_v$. 
Therefore there exists an immersion $F$ of a neighborhood of a point 
of the $uv$-plane into $N$ with $K\equiv L_0$, $R^{\perp} \equiv 0$ 
such that the second fundamental form of $F$ does not satisfy 
the linearly dependent condition. 
In the above arguments, 
even if we replace \eqref{e=+rho} by \eqref{e=-rho} 
and $\varepsilon =1$ by $\varepsilon =-1$, 
we obtain the same conclusion. 

Hence we obtain 

\begin{thm}\label{thm:notldL0tn} 
Let $F:M\longrightarrow N$ be a time-like and conformal immersion 
satisfying $K\equiv L_0$ and $R^{\perp} \equiv 0$. 
Then $F$ satisfies $k_{\pm} \not= 0, \pm 1$ and \eqref{X+Y-X-Y+} 
if and only if 
$\lambda$ with $\lambda_{uu} -\lambda_{vv} +L_0 e^{2\lambda} =0$, 
$\gamma$, $t_-$, $\rho$, $f_+$, $f_-$ satisfy 
$t_- \not= 0$, \eqref{f+f-f+f-}, 
either \eqref{e=+rho} or \eqref{e=-rho}, and \eqref{dgamma1_2tn}. 
In addition, 
for functions $\lambda$, $\gamma$, $t_-$, $\rho$, $f_+$, $f_-$ 
as above defined on a simply connected open set $O$ 
of the $uv$-plane, 
there exists a time-like and conformal immersion $F$ of $O$ into $N$ 
with $K\equiv L_0$, $R^{\perp} \equiv 0$ 
such that the second fundamental form does not satisfy 
the linearly dependent condition, 
which is unique up to an isometry of $N$. 
\end{thm} 

In the following, suppose $L_0 =0$. 
Then we can suppose $\lambda \equiv 1$. 
Then \eqref{dgamma1_2tn} is represented as 
\begin{equation} 
 \left[ 
 \begin{array}{c} 
  \gamma_u \\ 
  \gamma_v 
   \end{array} 
 \right] =
 \left[ 
 \begin{array}{c} 
  (t_- )_u \\ 
  (t_- )_v 
   \end{array} 
 \right] 
-\delta 
 \dfrac{\dfrac{\partial (f_+ , \rho )}{\partial (u, v)}}{
        \dfrac{\partial (f_+ , f_-  )}{\partial (u, v)}}  
 \left[ 
 \begin{array}{c} 
  (f_- )_u \\ 
  (f_- )_v 
   \end{array} 
 \right] , 
\label{dgamma2tn} 
\end{equation} 
Since $\gamma_{uv} =\gamma_{vu}$, we obtain 
\begin{equation*} 
d\left( \dfrac{\dfrac{\partial (f_+ , \rho )}{\partial (u, v)}}{
               \dfrac{\partial (f_+ , f_-  )}{\partial (u, v)}} 
 \right) 
 \wedge df_-  
=0. 
\label{gammauvtn} 
\end{equation*} 
Therefore there exists a function $\Tilde{\xi}$ of one variable 
satisfying 
\begin{equation}
 \dfrac{\partial (f_+ , \rho )}{\partial (u, v)}
=\Tilde{\xi} (f_- )
 \dfrac{\partial (f_+ , f_-  )}{\partial (u, v)} . 
\label{gammauv2tn}
\end{equation}
From \eqref{e=+rho}, $(f_+ )_{uv} =(f_+ )_{vu}$ and \eqref{gammauv2tn}, 
we obtain 
\begin{equation} 
 \left[ 
 \begin{array}{c} 
  \rho_u \\ 
  \rho_v 
   \end{array} 
 \right] 
=({\rm cosh}^2\,\rho                 )\mbox{\boldmath{$a$}} 
+({\rm cosh}\,\rho\,{\rm sinh}\,\rho )\mbox{\boldmath{$b$}}
+({\rm sinh}^2\,\rho                 )\mbox{\boldmath{$c$}} , 
\label{gammauv3tn}
\end{equation} 
where 
\begin{equation} 
\begin{split} 
  \mbox{\boldmath{$a$}} := 
& \dfrac{2}{C} 
    \left[ 
    \begin{array}{cc} 
     \delta (f_- )_u &  (f_- )_v \\ 
    -\delta (f_- )_v & -(f_- )_u 
      \end{array}
    \right] 
    \left[ 
    \begin{array}{c} 
                 (f_- )_{uv} \\ 
    -\Tilde{\xi} (f_- )(f_- )_u  (f_- )_v 
      \end{array}
    \right] , \\ 
  \mbox{\boldmath{$b$}} :=
& \dfrac{1}{C} 
  \left( 
    \left[ 
    \begin{array}{cc} 
     (f_- )_u &  \delta  (f_- )_v \\ 
    -(f_- )_v & -\delta (f_- )_u 
      \end{array}
    \right] 
    \left[ 
    \begin{array}{c} 
    -(f_- )_{uu} -(f_- )_{vv} \\ 
     \Tilde{\xi} (f_- )((f_- )^2_u +(f_- )^2_v )
      \end{array}
    \right] 
  \right. \\ 
& \quad \quad 
  \left.   
  +2\left[ 
    \begin{array}{cc} 
    -(f_- )_v & -\delta (f_- )_u \\ 
     (f_- )_u &  \delta (f_- )_v 
      \end{array}
    \right] 
    \left[ 
    \begin{array}{c} 
                 (f_- )_{uv} \\ 
    -\Tilde{\xi} (f_- )(f_- )_u  (f_- )_v 
      \end{array}
    \right] 
    \right) , \\ 
  \mbox{\boldmath{$c$}} := 
& \dfrac{1}{C} 
    \left[ 
    \begin{array}{cc} 
    -\delta (f_- )_v & -(f_- )_u \\ 
     \delta (f_- )_u &  (f_- )_v 
      \end{array}
    \right] 
    \left[ 
    \begin{array}{c} 
    -(f_- )_{uu} -(f_- )_{vv} \\ 
     \Tilde{\xi} (f_- )((f_- )^2_u +(f_- )^2_v 
      \end{array}
    \right] . 
\end{split}
\label{abctn}
\end{equation}
We can represent \eqref{gammauv3tn} as in \eqref{abc2} 
with $t:={\rm tanh}\,\rho$, \eqref{omega} and \eqref{abctn}. 

Let $\Tilde{\xi}$ be a function of one variable. 
Let $f_-$ be a function of two variables $u$, $v$ 
satisfying $(f_- )^2_u -(f_- )^2_v \not= 0$. 
Let $t$ be a solution of \eqref{abc2} valued in $(-1, 1)\setminus \{ 0\}$ 
satisfying 
\begin{equation} 
2\delta(f_- )_u (f_- )_v -((f_-)^2_u +(f_- )^2_v )t\not= 0. 
\label{ttn} 
\end{equation} 
Let $\rho$ be a function given by $t={\rm tanh}\,\rho$. 
Then there exists a function $f_+$ 
satisfying $A\not= 0$, \eqref{e=+rho} and \eqref{gammauv2tn}. 
Since $f_-$ and $t$ satisfy \eqref{ttn}, we obtain $A'\not= 0$. 
Therefore there exist functions $\gamma$, $t_-$ 
satisfying $t_- \not= 0$ and \eqref{dgamma2tn}. 
Then $\lambda \equiv 1$, $\gamma$, $t_-$, $\rho$, $f_+$, $f_-$ 
satisfy \eqref{f+f-f+f-} and \eqref{dgamma1_2tn}. 

From \eqref{e=-rho}, $(f_+ )_{uv} =(f_+ )_{vu}$ and \eqref{gammauv2tn}, 
we obtain 
\begin{equation} 
 \left[ 
 \begin{array}{c} 
  \rho_u \\ 
  \rho_v 
   \end{array} 
 \right] 
=({\rm cosh}^2\,\rho                 )\mbox{\boldmath{$c$}} 
+({\rm cosh}\,\rho\,{\rm sinh}\,\rho )\mbox{\boldmath{$b$}}
+({\rm sinh}^2\,\rho                 )\mbox{\boldmath{$a$}} , 
\label{gammauv3tn-}
\end{equation} 
where $\mbox{\boldmath{$a$}}$, 
      $\mbox{\boldmath{$b$}}$, 
      $\mbox{\boldmath{$c$}}$ are as in \eqref{abctn}. 
We can represent \eqref{gammauv3tn-} as in \eqref{abc2} 
with $t:={\rm tanh}\,\rho$ and 
\begin{equation}
\omega_0 =[du \ dv]\mbox{\boldmath{$c$}} , \quad 
\omega_1 =[du \ dv]\mbox{\boldmath{$b$}} , \quad 
\omega_2 =[du \ dv]\mbox{\boldmath{$a$}} 
\label{omega-} 
\end{equation} 
and \eqref{abctn}. 

Let $\Tilde{\xi}$ be a function of one variable. 
Let $f_-$ be a function of two variables $u$, $v$ 
satisfying $(f_- )^2_u -(f_- )^2_v \not= 0$. 
Let $t$ be a solution of \eqref{abc2} with \eqref{omega-} 
valued in $(-1, 1)\setminus \{ 0\}$ satisfying 
\begin{equation} 
2\delta(f_- )_u (f_- )_v t -((f_-)^2_u +(f_- )^2_v )\not= 0. 
\label{ttn-} 
\end{equation} 
Let $\rho$ be a function given by $t={\rm tanh}\,\rho$. 
Then there exists a function $f_+$ 
satisfying $A\not= 0$, \eqref{e=-rho} and \eqref{gammauv2tn}. 
Since $f_-$ and $t$ satisfy \eqref{ttn-}, we obtain $A'\not= 0$. 
Therefore there exist functions $\gamma$, $t_-$ 
satisfying $t_- \not= 0$ and \eqref{dgamma2tn}. 
Then $\lambda \equiv 1$, $\gamma$, $t_-$, $\rho$, $f_+$, $f_-$ 
satisfy \eqref{f+f-f+f-} and \eqref{dgamma1_2tn}. 

Therefore from Theorem~\ref{thm:notldL0tn}, we obtain 

\begin{thm}\label{thm:notldtn} 
Let $F:M\longrightarrow E^4_2$ be a time-like and conformal immersion 
satisfying $K\equiv 0$ and $R^{\perp} \equiv 0$. 
Then $F$ satisfies $k_{\pm} \not= 0, \pm 1$ and \eqref{X+Y-X-Y+} 
if and only if 
the functions $\gamma$, $t_-$, $f_-$, $t={\rm tanh}\,\rho$ 
satisfy $t_- \not= 0$, $(f_- )^2_u -(f_- )^2_v \not= 0$, \eqref{abc2} 
with either \eqref{omega} or \eqref{omega-}, 
\eqref{dgamma2tn} and \eqref{abctn} 
for a function $\Tilde{\xi}$ of one variable 
and either \eqref{ttn} or \eqref{ttn-}. 
In addition, 
for functions $f_-$, $\rho$ as above 
defined on a simply connected open set $O$ of the $uv$-plane, 
there exists a time-like and conformal immersion $F$ of $O$ into $E^4_2$ 
with $K\equiv 0$, $R^{\perp} \equiv 0$ 
such that the second fundamental form does not satisfy 
the linearly dependent condition. 
\end{thm} 

\section{Space-like surfaces with flat normal connection 
in Lorentzian 4-dimensional space forms}\label{sect:sL} 

\setcounter{equation}{0} 

\subsection{The Gauss-Codazzi-Ricci equations}\label{subsect:sLGCR} 

Let $N$ be a 4-dimensional Lorentzian space form 
with constant sectional curvature $L_0$. 
If $L_0 =0$, 
then we can suppose $N=E^4_1$; 
if $L_0 >0$, 
then we can suppose $N=\{ x\in E^5_1 \ | \ \langle x, x\rangle =1/L_0 \}$; 
if $L_0 <0$, 
then we can suppose $N=\{ x\in E^5_2 \ | \ \langle x, x\rangle =1/L_0 \}$, 
where $\langle \ , \ \rangle$ is the metric of $E^5_1$ or $E^5_2$.  
Let $h$ be the metric of $N$. 
Let $M$ be a Riemann surface 
and $F:M\longrightarrow N$ a space-like and conformal immersion 
of $M$ into $N$. 
Let $(u, v)$ be as in Section~\ref{sect:R}. 
The induced metric $g$ on $M$ by $F$ is represented 
as $g=e^{2\lambda} (du^2 +dv^2 )$ 
for a real-valued function $\lambda$. 
Let $\Tilde{\nabla}$ denote the Levi-Civita connection 
of $E^4_1$, $E^5_1$ or $E^5_2$ according to $L_0 =0$, $>0$ or $<0$. 
We set $T_1 :=dF(\partial_u )$, $T_2 :=dF(\partial_v )$. 
Let $N_1$, $N_2$ be normal vector fields of $F$ satisfying 
$$h(N_1 , N_1 )=-h(N_2 , N_2 )=e^{2\lambda} , \quad 
  h(N_1 , N_2 )=0.$$ 
Then we have 
\begin{equation} 
\begin{split} 
\Tilde{\nabla}_{T_1} (T_1 \ T_2 \ N_1 \ N_2 \ F) 
& =                  (T_1 \ T_2 \ N_1 \ N_2 \ F)S, \\ 
\Tilde{\nabla}_{T_2} (T_1 \ T_2 \ N_1 \ N_2 \ F) 
& =                  (T_1 \ T_2 \ N_1 \ N_2 \ F)T, 
\end{split} 
\label{dt1t2sL} 
\end{equation} 
where 
\begin{equation} 
\begin{split} 
S & =\left[ \begin{array}{ccccc} 
         \lambda_u         & \lambda_v & -\alpha_1  &  \beta_1   & 1 \\ 
        -\lambda_v         & \lambda_u & -\alpha_2  &  \beta_2   & 0 \\ 
         \alpha_1          & \alpha_2  &  \lambda_u &  \mu_1     & 0 \\ 
         \beta_1           & \beta_2   &  \mu_1     &  \lambda_u & 0 \\ 
         -L_0 e^{2\lambda} &  0        &   0        &   0        & 0 
              \end{array} 
     \right] , \\ 
T & =\left[ \begin{array}{ccccc} 
         \lambda_v & -\lambda_u        & -\alpha_2  &  \beta_2   & 0 \\ 
         \lambda_u &  \lambda_v        & -\alpha_3  &  \beta_3   & 1 \\ 
         \alpha_2  &  \alpha_3         &  \lambda_v &  \mu_2     & 0 \\ 
         \beta_2   &  \beta_3          &  \mu_2     &  \lambda_v & 0 \\ 
          0        & -L_0 e^{2\lambda} &   0        &   0        & 0 
              \end{array} 
     \right] , 
\end{split} 
\label{STsL} 
\end{equation} 
and $\alpha_k$, $\beta_k$ ($k=1, 2, 3$) and $\mu_l$ ($l=1, 2$) are 
real-valued functions. 
From \eqref{dt1t2sL}, we obtain $S_v -T_u =ST-TS$. 
This is equivalent to 
the system of the equations of Gauss, Codazzi and Ricci. 
The equation of Gauss is given by 
\begin{equation} 
  \lambda_{uu} +\lambda_{vv} +L_0 e^{2\lambda} 
=-\alpha_1 \alpha_3 +\beta_1 \beta_3 +\alpha^2_2 -\beta^2_2 . 
\label{gausssL}
\end{equation}
The equations of Codazzi are given by 
\begin{equation} 
\begin{split}
   (\alpha_1 )_v -(\alpha_2 )_u 
& = \alpha_2 \lambda_u +\alpha_3 \lambda_v 
   +\beta_2 \mu_1      -\beta_1  \mu_2  , \\ 
   (\alpha_2 )_v -(\alpha_3 )_u 
& =-\alpha_1 \lambda_u -\alpha_2 \lambda_v 
   +\beta_3  \mu_1     -\beta_2  \mu_2  , \\ 
   (\beta_1 )_v -(\beta_2 )_u 
& = \beta_2  \lambda_u +\beta_3  \lambda_v 
   +\alpha_2 \mu_1     -\alpha_1 \mu_2  , \\ 
   (\beta_2 )_v -(\beta_3 )_u 
& =-\beta_1  \lambda_u -\beta_2  \lambda_v 
   +\alpha_3 \mu_1     -\alpha_2 \mu_2 .  
\end{split} 
\label{codazzisL} 
\end{equation} 
The equation of Ricci is given by 
\begin{equation} 
  (\mu_1 )_v -(\mu_2 )_u = \alpha_1 \beta_2 -\alpha_2 \beta_1 
                          +\alpha_2 \beta_3 -\alpha_3 \beta_2 . 
\label{riccisL} 
\end{equation}

We have analogues of Theorems~\ref{thm:pnvftn} and \ref{thm:pnvftnll}. 

\begin{rem}\label{rem:K=L_0sL} 
Suppose that $F$ has zero mean curvature vector. 
Then the following are mutually equivalent (\cite{ando3}, \cite{ando6}): 
\begin{itemize} 
\item[{\rm (a)}]{$K\equiv L_0$;} 
\item[{\rm (b)}]{the holomorphic quartic differential $Q$ vanishes;} 
\item[{\rm (c)}]{the covariant derivatives of 
the two lifts of $F$ are zero or light-like, 
and perpendicular to each other;} 
\item[{\rm (d)}]{a light-like normal vector field of $F$ is parallel.} 
\end{itemize} 
The conformal Gauss maps of Willmore surfaces in 
$3$-dimensional Riemannian (respectively, Lorentzian) space forms  
with vanishing holomorphic quartic differential satisfy 
the above (a) $\sim$ (d) with $L_0 =1$ (respectively, $L_0 =-1$) in (a) 
(\cite{bryant2}, \cite{ando3}, \cite{ando6}). 
In particular, the conformal Gauss maps 
of minimal (respectively, maximal) surfaces 
in Euclidean $3$-space $E^3$ (respectively, $E^3_1$) 
satisfy these conditions with $L_0 =1$ (respectively, $L_0 =-1$) in (a), 
and we can refer to \cite{ando0} (respectively, \cite{AHHK}) 
for such surfaces.  
If $F$ satisfies (d), 
then the second fundamental form of $F$ satisfies 
the light-like linearly dependent condition 
and therefore by \eqref{riccisL}, 
the normal connection $\nabla^{\perp}$ is flat. 
There exist space-like surfaces in $N$ 
with zero mean curvature vector, 
$K\not= L_0$ and $R^{\perp} \equiv 0$, 
and by the analogue of Theorem~\ref{thm:pnvftn}, 
such surfaces have space-like or time-like parallel normal vector fields. 
\end{rem} 

\subsection{Space-like surfaces with \mbox{\boldmath{$R^{\perp} \equiv 0$}} 
and \mbox{\boldmath{$K\equiv L_0$}} 
which do not admit any parallel normal vector fields}\label{subsect:sLnpnvf} 

We have two analogues of Theorem~\ref{thm:ldnotpnvf} 
for space-like surfaces in $N$, 
according to the space-like linearly dependent condition  
and          the time-like  linearly dependent condition. 

In the following, 
we study a space-like and conformal immersion $F:M\longrightarrow N$ 
of $M$ into $N$ with $R^{\perp} \equiv 0$ and $K\equiv L_0$ 
such that the second fundamental form of $F$ does not satisfy 
the linearly dependent condition. 

We set 
\begin{equation*} 
\begin{array}{lcl} 
W:=\alpha_2 -i\beta_1  , & \ & X:=\alpha_2 +i\beta_3 , \\ 
Y:=\beta_2  -i\alpha_1 , & \ & Z:=\beta_2  +i\alpha_3 , 
\end{array} 
\label{WXYZsL} 
\end{equation*} 
where $i$ is the imaginary unit of complex numbers. 
Then we have 
\begin{equation} 
W +\overline{W} =X +\overline{X} , \quad 
Y +\overline{Y} =Z +\overline{Z} . 
\label{W=XY=ZsL} 
\end{equation} 
By \eqref{gausssL} and \eqref{riccisL}, we obtain 
\begin{equation} 
 WX-YZ=\lambda_{uu} +\lambda_{vv} +L_0 e^{2\lambda} 
      +i((\mu_1 )_v -(\mu_2 )_u ). 
\label{gaussriccisL} 
\end{equation} 
We can rewrite \eqref{codazzisL} into 
\begin{equation} 
\begin{split}
   Y_v +iX_u 
& =-iW(\lambda_u -i\mu_2 )-Z(\lambda_v +i\mu_1 ), \\ 
   W_v +iZ_u 
& =-iY(\lambda_u -i\mu_2 )-X(\lambda_v +i\mu_1 ). 
\end{split} 
\label{codazziWXYZsL} 
\end{equation} 

Suppose $R^{\perp} \equiv 0$ and $K\equiv L_0$. 
Then from \eqref{gaussriccisL}, we obtain 
\begin{equation} 
WX-YZ=0. 
\label{WXYZ0sL} 
\end{equation} 
We suppose that there exists a complex-valued function $k$ 
satisfying $k\not= 0$ and $(W, Z)=k(Y, X)$.  
In addition, we suppose $k\not= \pm 1$. 
We will prove 

\begin{lem}\label{lem:fXYsL} 
There exists a complex-valued function $f$ satisfying 
\begin{equation} 
X=\dfrac{if_v}{\sqrt{k^2 -1}}, \quad 
Y=\dfrac{f_u}{\sqrt{k^2 -1}}. 
\label{fsL} 
\end{equation} 
\end{lem} 

\vspace{3mm} 

\par\noindent 
\textit{Proof} \ 
Applying $(W, Z)=k(Y, X)$ to \eqref{codazziWXYZsL}, we obtain 
\begin{equation} 
\begin{split}
   Y_v  +i  X_u 
& =-k(iY(\lambda_u -i\gamma_v )+X(\lambda_v +i\gamma_u )), \\ 
 (kY)_v +i(kX)_u 
& =-iY(\lambda_u -i\gamma_v )-X(\lambda_v +i\gamma_u ). 
\end{split} 
\label{codazziWXYZsL2} 
\end{equation} 
From \eqref{codazziWXYZsL2}, we obtain 
\begin{equation} 
kk_v Y+ikk_u X +(k^2 -1)(Y_v +iX_u )=0. 
\label{codazziWXYZsL2_2} 
\end{equation} 
We can consider a square root $\sqrt{k^2 -1}$ of $k^2 -1$ 
on a neighborhood of each point, 
and therefore we can rewrite \eqref{codazziWXYZsL2_2} into 
\begin{equation*} 
\left( \sqrt{k^2 -1} Y\right)_v +i\left( \sqrt{k^2 -1} X\right)_u =0. 
\label{codazziWXYZsL2_2_2} 
\end{equation*} 
This means that 
there exists a complex-valued function $f$ satisfying \eqref{fsL}. 
Hence we obtain Lemma~\ref{lem:fXYsL}. 
\hfill 
$\square$ 

\vspace{3mm} 

We will prove 

\begin{lem}\label{lem:f+2f-2sL} 
The complex-valued function $f$ satisfies 
\begin{equation} 
{\rm Re}\,f_u\,{\rm Im}\,f_u +{\rm Re}\,f_v\,{\rm Im}\,f_v =0. 
\label{fufvsL} 
\end{equation} 
\end{lem} 

\vspace{3mm} 

\par\noindent 
\textit{Proof} \ 
Since $R^{\perp} \equiv 0$ is equivalent to 
\begin{equation} 
           W^2  +          X^2  +\overline{Y}^2 +\overline{Z}^2 
=\overline{W}^2 +\overline{X}^2 +          Y^2  +          Z^2, 
\label{X2Y2Z2W2sL} 
\end{equation} 
from $(W, Z)=k(Y, X)$, \eqref{fsL} and \eqref{X2Y2Z2W2sL}, 
we obtain \eqref{fufvsL}. 
\hfill 
$\square$ 

\vspace{3mm} 

\begin{rem} 
The relation \eqref{fufvsL} is equivalent to 
a condition that $f^2_u +f^2_v$ is real-valued. 
\end{rem} 

Noticing $k\not= \pm 1$ and 
referring to the proof of Lemma~\ref{lem:XYXY},  
we can prove 

\begin{lem}\label{lem:XYXYsL} 
If $X^2 \overline{Y}^2 -\overline{X}^2 Y^2 \not= 0$, 
then the second fundamental form of $F$ does not satisfy 
the linearly dependent condition. 
\end{lem} 

In the following, 
we suppose $X^2 \overline{Y}^2 -\overline{X}^2 Y^2 \not= 0$, 
which is equivalent to 
\begin{equation} 
{\rm Re}\,(f_u \overline{f}_v ) 
{\rm Im}\,(f_u \overline{f}_v )\not= 0. 
\label{fufv} 
\end{equation} 
We set 
\begin{equation} 
A :=2{\rm Re}\,(f_u \overline{f}_v ), \ 
A':=2{\rm Im}\,(f_u \overline{f}_v ), \ 
B :=|f_u |^2 -|f_v |^2 ,              \ 
C := f^2_u    +f^2_v . 
\label{AA'BCsL} 
\end{equation} 
Then $A$, $B$, $C$ satisfy $A^2 +B^2 =C^2$, and 
\eqref{fufv} is equivalent to $AA' \not= 0$. 
In particular, we have $A\not= 0$, which yields $C\not= 0$. 
We will prove 

\begin{lem}\label{lem:k+k-sL} 
The complex-valued function $k$ satisfies 
\begin{equation} 
Ak+iB\not= 0, \quad 
\overline{k} =\dfrac{iBk +A}{Ak+iB} . 
\label{overlinek} 
\end{equation} 
\end{lem} 

\vspace{3mm} 

\par\noindent 
\textit{Proof} \ 
Applying $(W, Z)=k(Y, X)$ to \eqref{W=XY=ZsL}, 
we obtain 
\begin{equation} 
k=\dfrac{-X\overline{X} -\overline{X}^2 +Y\overline{Y} +\overline{Y}^2}{
          X\overline{Y} -\overline{X} Y} . 
\label{ksL} 
\end{equation} 
Applying \eqref{fsL} to \eqref{ksL}, we obtain 
\begin{equation} 
Aik-B=C\sqrt{\dfrac{k^2 -1}{\overline{k}^2 -1}} . 
\label{ABCsL} 
\end{equation} 
Since $C\not= 0$, from \eqref{ABCsL}, we obtain \eqref{overlinek}. 
\hfill 
$\square$ 

\vspace{3mm} 

From \eqref{overlinek}, we obtain 
\begin{equation} 
\overline{k}^2 -1=-C^2 \dfrac{k^2 -1}{(Ak +iB)^2} 
\label{overlinek2-1} 
\end{equation} 
and 
\begin{equation} 
\overline{k}_a 
=i(AB_a -BA_a )\dfrac{k^2 -1}{(Ak +iB)^2} 
-  C^2 \dfrac{k_a}{(Ak +iB)^2} . 
\label{overlineka} 
\end{equation} 
for $a=u, v$. 

From the first equation in \eqref{codazziWXYZsL2}, 
we obtain 
\begin{equation} 
\begin{split} 
&  \left[ 
   \begin{array}{cc} 
              k            X & -i          k            Y \\ 
    \overline{k} \overline{X} & i\overline{k} \overline{Y} 
     \end{array} 
   \right] 
   \left[ 
   \begin{array}{c} 
    \gamma_u \\ 
    \gamma_v 
     \end{array} 
   \right] \\ 
& =\left[ 
   \begin{array}{c} 
   -          X_u  +i          Y_v \\ 
   -\overline{X}_u -i\overline{Y}_v 
     \end{array} 
   \right]  
  +\left[ 
   \begin{array}{cc} 
   -          k            Y  &  i          k            X \\ 
   -\overline{k} \overline{Y} & -i\overline{k} \overline{X} 
     \end{array} 
   \right] 
   \left[ 
   \begin{array}{c} 
    \lambda_u \\ 
    \lambda_v 
     \end{array} 
   \right] . 
\end{split} 
\label{gammasL} 
\end{equation} 
Then applying \eqref{fsL}, \eqref{overlinek2-1} and \eqref{overlineka} 
to \eqref{gammasL}, we obtain 
\begin{equation} 
\begin{split} 
    \left[ 
    \begin{array}{c} 
     \gamma_u \\ 
     \gamma_v 
      \end{array} 
    \right] =
&   \dfrac{1}{k^2 -1} 
    \left[ 
    \begin{array}{c} 
     k_u \\ 
     k_v 
      \end{array} 
    \right] 
  +i\dfrac{A^2}{C^2} 
   \dfrac{\dfrac{\partial (\overline{f}, B/A )}{\partial (u, v)}}{
          \dfrac{\partial (f, \overline{f} )}{\partial (u, v)}} 
   \left[ 
   \begin{array}{c} 
    f_u \\ 
    f_v 
     \end{array} 
   \right] \\ 
& -\dfrac{1}{A'} 
   \left[ 
   \begin{array}{cc} 
    2|f_u |^2 &   A \\ 
      A       & 2|f_v |^2 
     \end{array} 
   \right] 
   \left[ 
   \begin{array}{c} 
    \lambda_u \\ 
    \lambda_v 
     \end{array} 
   \right] . 
\end{split} 
\label{dgammasL} 
\end{equation} 
By \eqref{fufvsL}, 
there exists a real-valued function $\psi$ satisfying 
\begin{equation} 
 \left[ \begin{array}{c} 
         \overline{f}_u \\ 
         \overline{f}_v 
          \end{array} 
 \right] 
=\left[ \begin{array}{cc} 
         \cos \psi & -\sin \psi \\ 
         \sin \psi &  \cos \psi 
          \end{array} 
 \right] 
 \left[ \begin{array}{c} 
         f_v \\ 
         f_u 
          \end{array} 
 \right] . 
\label{psisL} 
\end{equation} 
Then we have $A=C\cos \psi$ and $B=-C\sin \psi$. 
Therefore $B/A =-\tan \psi$. 
Then \eqref{dgammasL} is represented as 
\begin{equation} 
\begin{split} 
    \left[ 
    \begin{array}{c} 
     \gamma_u \\ 
     \gamma_v 
      \end{array} 
    \right] = 
&   \dfrac{1}{k^2 -1} 
    \left[ 
    \begin{array}{c} 
     k_u \\ 
     k_v 
      \end{array} 
    \right] 
  -i\dfrac{\dfrac{\partial (\overline{f}, \psi )}{\partial (u, v)}}{
           \dfrac{\partial (f, \overline{f} )}{\partial (u, v)}} 
    \left[ 
    \begin{array}{c} 
     f_u \\ 
     f_v 
      \end{array} 
    \right] \\ 
&  -\dfrac{1}{A'} 
    \left[ 
    \begin{array}{cc} 
     2|f_u |^2 &   A \\ 
       A       & 2|f_v |^2 
      \end{array} 
    \right] 
    \left[ 
    \begin{array}{c} 
     \lambda_u \\ 
     \lambda_v 
      \end{array} 
    \right] . 
\end{split} 
\label{dgamma1_2sL} 
\end{equation} 

Let $\lambda$ be a function 
satisfying $\lambda_{uu} +\lambda_{vv} +L_0 e^{2\lambda} =0$. 
Let $\gamma$, $\psi$ be real-valued functions 
and $k$, $f$ complex-valued functions satisfying $k\not= 0, \pm 1$, 
\eqref{fufvsL}, \eqref{fufv}, \eqref{overlinek}, \eqref{psisL} 
and \eqref{dgamma1_2sL}. 
Then the above functions satisfy \eqref{dgammasL}. 
We define $W$, $X$, $Y$, $Z$ by 
$(W, Z)=k(Y, X)$ and \eqref{fsL}. 
Then these functions satisfy \eqref{W=XY=ZsL} and \eqref{gaussriccisL} 
with $\mu_1 =\gamma_u$ and $\mu_2 =\gamma_v$. 
In addition, they satisfy \eqref{codazziWXYZsL} 
and $X^2 \overline{Y}^2 -\overline{X}^2 Y^2 \not= 0$. 
Noticing \eqref{W=XY=ZsL}, 
we define functions $\alpha_j$, $\beta_j$ ($j=1, 2, 3$) by 
\begin{equation*} 
\alpha_1 :=\dfrac{i}{2}  (Y -\overline{Y} ), \quad 
\alpha_2 :=\dfrac{1}{2}  (X +\overline{X} ), \quad 
\alpha_3 :=\dfrac{1}{2i} (Z -\overline{Z} ) 
\label{alpha123sL} 
\end{equation*} 
and 
\begin{equation*} 
\beta_1 :=\dfrac{i}{2}  (W -\overline{W} ), \quad 
\beta_2 :=\dfrac{1}{2}  (Y +\overline{Y} ), \quad 
\beta_3 :=\dfrac{1}{2i} (X -\overline{X} ). 
\label{beta123sL} 
\end{equation*} 
Then $\alpha_j$, $\beta_j$ ($j=1, 2, 3$) satisfy $S_v -T_u =ST-TS$ 
for $S$, $T$ in \eqref{STsL} with $\mu_1 =\gamma_u$ and $\mu_2 =\gamma_v$. 
Therefore there exists an immersion $F$ of a neighborhood of a point 
of the $uv$-plane into $N$ with $K\equiv L_0$, $R^{\perp} \equiv 0$ 
such that the second fundamental form of $F$ does not satisfy 
the linearly dependent condition. 

Hence we obtain 

\begin{thm}\label{thm:notldL0sL} 
Let $F:M\longrightarrow N$ be a space-like and conformal immersion 
satisfying $K\equiv L_0$, $R^{\perp} \equiv 0$ and $k\not= 0, \pm 1$. 
Then $F$ satisfies $X^2 \overline{Y}^2 -\overline{X}^2 Y^2 \not= 0$ 
if and only if 
$\lambda$ with $\lambda_{uu} +\lambda_{vv} +L_0 e^{2\lambda} =0$, 
$\gamma$, $k$, $\psi$, $f$ satisfy 
\eqref{fufvsL}, \eqref{fufv}, \eqref{overlinek}, \eqref{psisL} 
and \eqref{dgamma1_2sL}. 
In addition, 
for functions $\lambda$, $\gamma$, $k$, $\psi$, $f$ 
as above defined on a simply connected open set $O$ 
of the $uv$-plane, 
there exists a space-like and conformal immersion $F$ of $O$ into $N$ 
with $K\equiv L_0$, $R^{\perp} \equiv 0$ 
such that the second fundamental form does not satisfy 
the linearly dependent condition, 
which is unique up to an isometry of $N$. 
\end{thm} 

In the following, suppose $L_0 =0$. 
Then we can suppose $\lambda \equiv 1$. 
Then \eqref{dgamma1_2sL} is represented as 
\begin{equation} 
   \left[ 
   \begin{array}{c} 
    \gamma_u \\ 
    \gamma_v 
     \end{array} 
   \right] = 
   \dfrac{1}{k^2 -1} 
   \left[ 
   \begin{array}{c} 
    k_u \\ 
    k_v 
     \end{array} 
   \right] 
 -i\Phi 
   \left[ 
   \begin{array}{c} 
    f_u \\ 
    f_v 
     \end{array} 
   \right] , \quad 
\Phi :=\dfrac{\dfrac{\partial (\overline{f}, \psi )}{\partial (u, v)}}{
              \dfrac{\partial (f, \overline{f} )}{\partial (u, v)}} . 
\label{dgamma2sL} 
\end{equation} 
Computing the first and second terms of the right side of \eqref{dgamma2sL} 
by \eqref{overlinek2-1}, \eqref{overlineka} and \eqref{psisL}, 
we see that 
the imaginary part of the right side of \eqref{dgamma2sL} vanishes: 
\begin{equation} 
 {\rm Re}\,(\Phi f_a )
={\rm Im}\,\dfrac{k_a}{k^2 -1} 
=-\dfrac{1}{2} \psi_a 
\label{psia} 
\end{equation} 
for $a=u, v$. 
Since $\gamma_{uv} =\gamma_{vu}$, we obtain $d\Phi \wedge df=0$. 
Therefore there exist real-valued functions $r$, $\theta$ satisfying 
\begin{equation} 
(\Phi_u , \Phi_v )=re^{i\theta} (f_u , f_v ). 
\label{Phi} 
\end{equation} 
Since ${\rm Im}\,(f_u \overline{f}_v )\not= 0$, 
\begin{equation} 
s:= {\rm Re}\,f, \quad 
t:= {\rm Im}\,f 
\label{st} 
\end{equation} 
form local coordinates of $M$. 
Therefore, noticing that \eqref{Phi} can be rewritten into 
\begin{equation*} 
  \left[ 
  \begin{array}{cc} 
   {\rm Re}\,\Phi_u & {\rm Im}\,\Phi_u \\ 
   {\rm Re}\,\Phi_v & {\rm Im}\,\Phi_v 
    \end{array} 
  \right]  
=r\left[ 
  \begin{array}{cc} 
   s_u & t_u \\ 
   s_v & t_v 
    \end{array} 
  \right] 
  \left[ 
  \begin{array}{cc} 
   \cos \theta & \sin \theta \\ 
  -\sin \theta & \cos \theta  
    \end{array} 
  \right] , 
\end{equation*} 
we see by the Cauchy-Riemann equation that 
there exists a holomorphic function $\Tilde{\xi}$ 
of $s+it$ satisfying $\Phi =\Tilde{\xi} (f)$, i.e., 
\begin{equation} 
 \dfrac{\partial (\overline{f}, \psi )}{\partial (u, v)} 
=\Tilde{\xi} (f) 
 \dfrac{\partial (f, \overline{f} )}{\partial (u, v)} . 
 \label{xif} 
\end{equation}  
By \eqref{psia}, we have 
\begin{equation} 
\psi =-2\,{\rm Re}\,(\Tilde{\Xi} (s+it)), 
\label{psiXi} 
\end{equation} 
where $\Tilde{\Xi}$ is a primitive function of $\Tilde{\xi}$. 
We can rewrite \eqref{psisL} into 
\begin{equation} 
 \left[ 
 \begin{array}{cc} 
   s_u & -t_u \\ 
   s_v & -t_v 
   \end{array} 
 \right] 
=\left[ 
 \begin{array}{cc} 
   \cos \psi & -\sin \psi \\ 
   \sin \psi &  \cos \psi  
   \end{array} 
 \right] 
 \left[ 
  \begin{array}{cc} 
   s_v & t_v \\ 
   s_u & t_u 
    \end{array} 
  \right] . 
\label{psisL2} 
\end{equation} 
Noticing \eqref{fufvsL}, we rewrite \eqref{psisL2} into 
\begin{equation} 
  \left[ 
  \begin{array}{cc} 
    \cos \psi & -\sin \psi \\ 
    \sin \psi &  \cos \psi  
    \end{array} 
  \right] 
=-\dfrac{1}{s_u t_v -s_v t_u} 
  \left[ 
  \begin{array}{cc} 
    s_u t_u -s_v t_v & -s_u t_v -s_v t_u \\ 
    s_u t_v +s_v t_u &  s_u t_u -s_v t_v \\ 
    \end{array} 
  \right] . 
\label{psist} 
\end{equation} 
In addition, \eqref{psist} is equivalent to the following system: 
\begin{equation} 
\begin{split} 
(\cos \psi )s_u +(\sin \psi -1)s_v & =0, \\ 
(\cos \psi )t_u +(\sin \psi +1)t_v & =0. 
\end{split} 
\label{psist2} 
\end{equation} 

Let $\Tilde{\Xi}$ be a holomorphic function of one complex variable. 
Then noticing that $\psi$ in \eqref{psiXi} depends analytically on $s$, $t$, 
we obtain a unique analytic solution $(s, t)$ 
of the initial value problem for the system \eqref{psist2} 
(see \cite[Chapter 1]{CH} for the Cauchy-Kowalevski theorem). 
Let $f$ be a complex-valued function of two variables $u$, $v$ 
satisfying \eqref{fufvsL}, \eqref{fufv} and \eqref{psist2} 
with \eqref{st} and \eqref{psiXi}. 
Then $AA' \not= 0$, and $C=f^2_u +f^2_v$ is real-valued and nonzero. 
In addition, $f$ satisfies \eqref{xif}. 
Let $k$ be a complex-valued function 
satisfying $k\not= \pm 1$ and \eqref{overlinek}. 
Then $k$ is nonzero, 
and it satisfies \eqref{psia}. 
Therefore there exists a real-valued function $\gamma$ 
satisfying \eqref{dgamma2sL}. 
Then $\lambda \equiv 1$, $\gamma$, $k$, $\psi$, $f$ are functions 
satisfying \eqref{dgamma1_2sL}. 

Therefore from Theorem~\ref{thm:notldL0sL}, we obtain 

\begin{thm}\label{thm:notldsL} 
Let $F:M\longrightarrow E^4_1$ be a space-like and conformal immersion 
satisfying $K\equiv 0$, $R^{\perp} \equiv 0$ and $k\not= 0, \pm 1$. 
Then $F$ satisfies $X^2 \overline{Y}^2 -\overline{X}^2 Y^2 \not= 0$ 
if and only if 
the complex-valued function $f$ satisfies \eqref{fufvsL}, \eqref{fufv} and 
\eqref{psist2} with \eqref{st} and \eqref{psiXi} 
for a holomorphic function $\Tilde{\Xi}$. 
In addition, 
for functions $f$ and $\psi$ as above 
defined on a simply connected open set $O$ of the $uv$-plane, 
there exists a space-like and conformal immersion $F$ of $O$ into $E^4_1$ 
with $K\equiv 0$, $R^{\perp} \equiv 0$ 
such that the second fundamental form does not satisfy 
the linearly dependent condition. 
\end{thm} 

\section{Time-like surfaces with flat normal connection 
in Lorentzian 4-dimensional space forms}\label{sect:tL} 

\setcounter{equation}{0} 

\subsection{The Gauss-Codazzi-Ricci equations}\label{subsect:tLGCR} 

Let $N$, $h$ be as in Section~\ref{sect:sL}. 
Let $M$ be a Lorentz surface 
and $F:M\longrightarrow N$ a time-like and conformal immersion 
of $M$ into $N$. 
Let $(u, v)$ be as in Section~\ref{sect:tn}. 
Then the induced metric $g$ on $M$ by $F$ is represented 
as $g=e^{2\lambda} (du^2 -dv^2 )$ 
for a real-valued function $\lambda$. 
Let $\Tilde{\nabla}$, $T_1$, $T_2$ be as in Section~\ref{sect:sL}. 
Let $N_1$, $N_2$ be normal vector fields of $F$ satisfying 
$$h(N_1 , N_1 )=h(N_2 , N_2 )=e^{2\lambda} , \quad 
  h(N_1 , N_2 )=0.$$ 
Then we have 
\begin{equation} 
\begin{split} 
\Tilde{\nabla}_{T_1} (T_1 \ T_2 \ N_1 \ N_2 \ F) 
& =                  (T_1 \ T_2 \ N_1 \ N_2 \ F)S, \\ 
\Tilde{\nabla}_{T_2} (T_1 \ T_2 \ N_1 \ N_2 \ F) 
& =                  (T_1 \ T_2 \ N_1 \ N_2 \ F)T, 
\end{split} 
\label{dt1t2tL} 
\end{equation} 
where 
\begin{equation} 
\begin{split} 
S & =\left[ \begin{array}{ccccc} 
       \lambda_u         & \lambda_v & -\alpha_1  & -\beta_1   & 1 \\ 
       \lambda_v         & \lambda_u &  \alpha_2  &  \beta_2   & 0 \\ 
       \alpha_1          & \alpha_2  &  \lambda_u & -\mu_1     & 0 \\ 
       \beta_1           & \beta_2   &  \mu_1     &  \lambda_u & 0 \\ 
       -L_0 e^{2\lambda} &  0        &   0        &   0        & 0 
              \end{array} 
     \right] , \\ 
T & =\left[ \begin{array}{ccccc} 
       \lambda_v &  \lambda_u        & -\alpha_2  & -\beta_2   & 0 \\ 
       \lambda_u &  \lambda_v        &  \alpha_3  &  \beta_3   & 1 \\ 
       \alpha_2  &  \alpha_3         &  \lambda_v & -\mu_2     & 0 \\ 
       \beta_2   &  \beta_3          &  \mu_2     &  \lambda_v & 0 \\ 
        0        &  L_0 e^{2\lambda} &   0        &   0        & 0 
              \end{array} 
    \right] , 
\end{split} 
\label{STtL} 
\end{equation} 
and $\alpha_k$, $\beta_k$ ($k=1, 2, 3$) and $\mu_l$ ($l=1, 2$) are 
real-valued functions. 
From \eqref{dt1t2tL}, we obtain $S_v -T_u =ST-TS$. 
This is equivalent to 
the system of the equations of Gauss, Codazzi and Ricci. 
The equation of Gauss is given by 
\begin{equation} 
  \lambda_{uu} -\lambda_{vv} +L_0 e^{2\lambda} 
= \alpha_1 \alpha_3 +\beta_1 \beta_3 -\alpha^2_2 -\beta^2_2 . 
\label{gausstL}
\end{equation}
The equations of Codazzi are given by 
\begin{equation} 
\begin{split}
   (\alpha_1 )_v -(\alpha_2 )_u 
& = \alpha_2 \lambda_u -\alpha_3 \lambda_v 
   -\beta_2 \mu_1      +\beta_1  \mu_2  , \\ 
   (\alpha_2 )_v -(\alpha_3 )_u 
& = \alpha_1 \lambda_u -\alpha_2 \lambda_v 
   -\beta_3  \mu_1     +\beta_2  \mu_2  , \\ 
   (\beta_1 )_v -(\beta_2 )_u 
& = \beta_2  \lambda_u -\beta_3  \lambda_v 
   +\alpha_2 \mu_1     -\alpha_1 \mu_2  , \\ 
   (\beta_2 )_v -(\beta_3 )_u 
& = \beta_1  \lambda_u -\beta_2  \lambda_v 
   +\alpha_3 \mu_1     -\alpha_2 \mu_2 .  
\end{split} 
\label{codazzitL} 
\end{equation} 
The equation of Ricci is given by 
\begin{equation} 
  (\mu_1 )_v -(\mu_2 )_u = \alpha_1 \beta_2 -\alpha_2 \beta_1 
                          -\alpha_2 \beta_3 +\alpha_3 \beta_2 . 
\label{riccitL} 
\end{equation}

Theorem~\ref{thm:pnvf} holds 
for a time-like and conformal immersion $F:M\longrightarrow N$. 

\begin{rem} 
Suppose that $F$ has zero mean curvature vector. 
Then referring to \cite{ando6}, 
we can show that the following are mutually equivalent: 
\begin{itemize} 
\item[{\rm (a)}]{$K\equiv L_0$;} 
\item[{\rm (b)}]{the paraholomorphic quartic differential $Q$ is 
zero or null;} 
\item[{\rm (c)}]{the covariant derivatives of 
the two lifts of $F$ are zero or light-like, 
and perpendicular to each other;} 
\item[{\rm (d)}]{the shape operator of any normal vector field is 
zero or light-like.} 
\end{itemize} 
See \cite{ando9} for construction of time-like surfaces in $N$ 
with zero mean curvature vector and $K\equiv L_0$. 
If $F$ satisfies (d), 
then the second fundamental form of $F$ satisfies 
the linearly dependent condition 
and therefore by \eqref{riccitL}, 
the normal connection $\nabla^{\perp}$ is flat. 
There exist time-like surfaces in $N$ 
with zero mean curvature vector, 
$K\not= L_0$ and $R^{\perp} \equiv 0$, 
and by the analogue of Theorem~\ref{thm:pnvftn}, 
such surfaces have space-like parallel normal vector fields. 
\end{rem} 

\subsection{Time-like surfaces with \mbox{\boldmath{$R^{\perp} \equiv 0$}} 
and \mbox{\boldmath{$K\equiv L_0$}} 
which do not admit any parallel normal vector fields}\label{subsect:tLnpnvf} 

We have an analogue of Theorem~\ref{thm:ldnotpnvf} 
for time-like surfaces in $N$. 

In the following, 
we study a time-like and conformal immersion $F:M\longrightarrow N$ 
of $M$ into $N$ with $R^{\perp} \equiv 0$ and $K\equiv L_0$ 
such that the second fundamental form of $F$ does not satisfy 
the linearly dependent condition. 

We set 
\begin{equation*} 
\begin{array}{lcl} 
W:=\alpha_2 +i\beta_1  , & \ & X:=\alpha_2 +i\beta_3 , \\ 
Y:=\beta_2  -i\alpha_1 , & \ & Z:=\beta_2  -i\alpha_3 . 
\end{array} 
\label{WXYZtL} 
\end{equation*} 
Then we have \eqref{W=XY=ZsL}.  
By \eqref{gausstL} and \eqref{riccitL}, we obtain 
\begin{equation} 
-WX-YZ=\lambda_{uu} -\lambda_{vv} +L_0 e^{2\lambda} 
      +i((\mu_1 )_v -(\mu_2 )_u ). 
\label{gaussriccitL} 
\end{equation} 
We can rewrite \eqref{codazzitL} into 
\begin{equation} 
\begin{split}
   Y_v +iX_u 
& =-iW(\lambda_u -i\mu_2 )-Z(\lambda_v -i\mu_1 ), \\ 
   W_v -iZ_u 
& = iY(\lambda_u -i\mu_2 )-X(\lambda_v -i\mu_1 ). 
\end{split} 
\label{codazziWXYZtL} 
\end{equation} 

Suppose $R^{\perp} \equiv 0$ and $K\equiv L_0$. 
Then from \eqref{gaussriccitL}, we obtain 
\begin{equation} 
WX+YZ=0. 
\label{WXYZ0tL} 
\end{equation} 
We suppose that there exists a complex-valued function $k$ 
satisfying $k\not= 0$ and $(W, Z)=k(-Y, X)$. 
In addition, we suppose $k\not= \pm i$. 
We will prove 

\begin{lem}\label{lem:fXYtL} 
There exists a complex-valued function $f$ 
satisfying 
\begin{equation} 
X=\dfrac{if_v}{\sqrt{k^2 +1}}, \quad 
Y=\dfrac{f_u}{\sqrt{k^2 +1}}. 
\label{ftL} 
\end{equation} 
\end{lem} 

\vspace{3mm} 

\par\noindent 
\textit{Proof} \ 
Applying $(W, Z)=k(-Y, X)$ to \eqref{codazziWXYZtL}, we obtain 
\begin{equation} 
\begin{split}
   Y_v  +i  X_u 
& = k(iY(\lambda_u -i\gamma_v )-X(\lambda_v -i\gamma_u )), \\ 
 (kY)_v +i(kX)_u 
& =-  iY(\lambda_u -i\gamma_v )+X(\lambda_v -i\gamma_u ). 
\end{split} 
\label{codazziWXYZtL2} 
\end{equation} 
From \eqref{codazziWXYZtL2}, we obtain 
\begin{equation} 
kk_v Y+ikk_u X+(k^2 +1)(Y_v +iX_u )=0. 
\label{codazziWXYZtL2_2} 
\end{equation} 
Then we can consider a square root $\sqrt{k^2 +1}$ of $k^2 +1$ 
on a neighborhood of each point, 
and we can rewrite \eqref{codazziWXYZtL2_2} into 
\begin{equation*} 
\left( \sqrt{k^2 +1} Y\right)_v +i\left( \sqrt{k^2 +1} X\right)_u =0. 
\label{codazziWXYZtL2_2_2} 
\end{equation*} 
This means that 
there exists a complex-valued function $f$ satisfying \eqref{ftL}. 
Hence we obtain Lemma~\ref{lem:fXYtL}. 
\hfill 
$\square$ 

\vspace{3mm} 

We will prove 

\begin{lem}\label{lem:f+2f-2tL} 
The complex-valued function $f$ satisfies 
\begin{equation} 
{\rm Re}\,f_u\,{\rm Im}\,f_u = {\rm Re}\,f_v\,{\rm Im}\,f_v . 
\label{fufvtL} 
\end{equation} 
\end{lem} 

\vspace{3mm} 

\par\noindent 
\textit{Proof} \ 
Since $R^{\perp} \equiv 0$ is equivalent to 
\begin{equation} 
 \overline{W}^2 +\overline{X}^2 +\overline{Y}^2 +\overline{Z}^2 
=          W^2  +          X^2  +          Y^2  +          Z^2, 
\label{X2Y2Z2W2tL} 
\end{equation} 
from $(W, Z)=k(-Y, X)$, \eqref{ftL} and \eqref{X2Y2Z2W2tL}, 
we obtain \eqref{fufvtL}. 
\hfill 
$\square$ 

\vspace{3mm} 

\begin{rem} 
The relation \eqref{fufvtL} is equivalent to 
a condition that $f^2_u -f^2_v$ is real-valued. 
\end{rem} 

Noticing $k\not= \pm i$ and 
referring to the proof of Lemma~\ref{lem:XYXY}, 
we obtain Lemma~\ref{lem:XYXYsL} for the present situation. 

In the following, 
we suppose $X^2 \overline{Y}^2 -\overline{X}^2 Y^2 \not= 0$, 
which is equivalent to \eqref{fufv}. 
Let $A$, $A'$ be as in \eqref{AA'BCsL} and set 
\begin{equation*} 
B :=|f_u |^2 +|f_v |^2 , \quad 
C := f^2_u -f^2_v . 
\end{equation*} 
Then $A$, $B$, $C$ satisfy $B^2 -A^2 =C^2$, 
and \eqref{fufv} is equivalent to $AA' \not= 0$. 
In particular, we have $A' \not= 0$, which yields $C\not= 0$. 
We will prove 

\begin{lem}\label{lem:k+k-tL} 
The complex-valued function $k$ satisfies 
\begin{equation} 
Ak+iB\not= 0, \quad 
\overline{k} =\dfrac{iBk -A}{Ak+iB} . 
\label{overlinektL} 
\end{equation} 
\end{lem} 

\vspace{3mm} 

\par\noindent 
\textit{Proof} \ 
Applying $(W, Z)=k(-Y, X)$ to \eqref{W=XY=ZsL}, 
we obtain 
\begin{equation} 
k=\dfrac{X\overline{X} +\overline{X}^2 +Y\overline{Y} +\overline{Y}^2}{
         X\overline{Y} -\overline{X} Y} . 
\label{ktL} 
\end{equation} 
Applying \eqref{ftL} to \eqref{ktL}, we obtain 
\begin{equation} 
Aik-B=C\sqrt{\dfrac{k^2 +1}{\overline{k}^2 +1}} . 
\label{ABCtL} 
\end{equation} 
Since $C\not= 0$, from \eqref{ABCtL}, we obtain \eqref{overlinektL}. 
\hfill 
$\square$ 

\vspace{3mm} 

From \eqref{overlinektL}, we obtain 
\begin{equation} 
\overline{k}^2 +1=-C^2 \dfrac{k^2 +1}{(Ak +iB)^2} 
\label{overlinek2+1} 
\end{equation} 
and 
\begin{equation} 
\overline{k}_a 
=i(AB_a -BA_a )\dfrac{k^2 +1}{(Ak +iB)^2} 
-  C^2 \dfrac{k_a}{(Ak +iB)^2} . 
\label{overlinekatL} 
\end{equation} 
for $a=u, v$. 

From the first equation in \eqref{codazziWXYZtL2}, 
we obtain 
\begin{equation} 
\begin{split} 
&  \left[ 
   \begin{array}{cc} 
              k            X & -i          k            Y \\ 
    \overline{k} \overline{X} & i\overline{k} \overline{Y} 
     \end{array} 
   \right] 
   \left[ 
   \begin{array}{c} 
    \gamma_u \\ 
    \gamma_v 
     \end{array} 
   \right] \\ 
& =\left[ 
   \begin{array}{c} 
              X_u  -i          Y_v \\ 
    \overline{X}_u +i\overline{Y}_v 
     \end{array} 
   \right]  
  +\left[ 
   \begin{array}{cc} 
   -          k            Y  & -i          k            X \\ 
   -\overline{k} \overline{Y} &  i\overline{k} \overline{X} 
     \end{array} 
   \right] 
   \left[ 
   \begin{array}{c} 
    \lambda_u \\ 
    \lambda_v 
     \end{array} 
   \right] . 
\end{split} 
\label{gammatL} 
\end{equation} 
Then applying \eqref{ftL}, \eqref{overlinek2+1} and \eqref{overlinekatL} 
to \eqref{gammatL}, we obtain 
\begin{equation} 
\begin{split} 
    \left[ 
    \begin{array}{c} 
     \gamma_u \\ 
     \gamma_v 
      \end{array} 
    \right] =
&  -\dfrac{1}{k^2 +1} 
    \left[ 
    \begin{array}{c} 
     k_u \\ 
     k_v 
      \end{array} 
    \right] 
  -i\dfrac{A^2}{C^2} 
   \dfrac{\dfrac{\partial (\overline{f}, B/A )}{\partial (u, v)}}{
          \dfrac{\partial (f, \overline{f} )}{\partial (u, v)}} 
   \left[ 
   \begin{array}{c} 
    f_u \\ 
    f_v 
     \end{array} 
   \right] \\ 
& -\dfrac{1}{A'} 
   \left[ 
   \begin{array}{cc} 
    2|f_u |^2 & -  A \\ 
      A       & -2|f_v |^2 
     \end{array} 
   \right] 
   \left[ 
   \begin{array}{c} 
    \lambda_u \\ 
    \lambda_v 
     \end{array} 
   \right] . 
\end{split} 
\label{dgammatL} 
\end{equation} 
Noticing \eqref{fufvtL}, 
we see that there exists a real-valued function $\rho$ satisfying 
\begin{equation} 
 \left[ \begin{array}{c} 
         \overline{f}_u \\ 
         \overline{f}_v 
          \end{array} 
 \right] 
=\left[ \begin{array}{cc} 
         \delta\,{\rm cosh}\,\rho &         {\rm sinh}\,\rho \\ 
                 {\rm sinh}\,\rho & \delta\,{\rm cosh}\,\rho 
          \end{array} 
 \right] 
 \left[ \begin{array}{c} 
         f_v \\ 
        -f_u 
          \end{array} 
 \right] 
\label{rhotL} 
\end{equation} 
for $\delta \in \{ 1, -1\}$. 
Then we have $A=-\delta C\,{\rm cosh}\,\rho$ 
and          $B=-       C\,{\rm sinh}\,\rho$. 
Therefore  $B/A=   \delta\,{\rm tanh}\,\rho$. 
Then \eqref{dgammatL} is represented as 
\begin{equation} 
\begin{split} 
    \left[ 
    \begin{array}{c} 
     \gamma_u \\ 
     \gamma_v 
      \end{array} 
    \right] = 
&  -\dfrac{1}{k^2 +1} 
    \left[ 
    \begin{array}{c} 
     k_u \\ 
     k_v 
      \end{array} 
    \right] 
  -i\delta 
    \dfrac{\dfrac{\partial (\overline{f}, \rho )}{\partial (u, v)}}{
           \dfrac{\partial (f, \overline{f} )}{\partial (u, v)}} 
    \left[ 
    \begin{array}{c} 
     f_u \\ 
     f_v 
      \end{array} 
    \right] \\ 
&  -\dfrac{1}{A'} 
    \left[ 
    \begin{array}{cc} 
     2|f_u |^2 & -  A \\ 
       A       & -2|f_v |^2 
      \end{array} 
    \right] 
    \left[ 
    \begin{array}{c} 
     \lambda_u \\ 
     \lambda_v 
      \end{array} 
    \right] . 
\end{split} 
\label{dgamma1_2tL} 
\end{equation} 

Let $\lambda$ be a function 
satisfying $\lambda_{uu} -\lambda_{vv} +L_0 e^{2\lambda} =0$. 
Let $\gamma$, $\rho$ be real-valued functions 
and $k$, $f$ complex-valued functions satisfying $k\not= 0, \pm i$, 
\eqref{fufv}, \eqref{fufvtL}, \eqref{overlinektL}, \eqref{rhotL} 
and \eqref{dgamma1_2tL}. 
Then the above functions satisfy \eqref{dgammatL}. 
We define $W$, $X$, $Y$, $Z$ by 
$(W, Z)=k(-Y, X)$ and \eqref{ftL}. 
Then these functions satisfy \eqref{W=XY=ZsL} and \eqref{gaussriccitL} 
with $\mu_1 =\gamma_u$ and $\mu_2 =\gamma_v$. 
In addition, they satisfy \eqref{codazziWXYZtL} 
and $X^2 \overline{Y}^2 -\overline{X}^2 Y^2 \not= 0$. 
Noticing \eqref{W=XY=ZsL}, 
we define functions $\alpha_k$, $\beta_k$ ($k=1, 2, 3$) by 
\begin{equation*} 
\alpha_1 :=\dfrac{i}{2}  (Y -\overline{Y} ), \quad 
\alpha_2 :=\dfrac{1}{2}  (X +\overline{X} ), \quad 
\alpha_3 :=\dfrac{i}{2} (Z -\overline{Z} ) 
\label{alpha123tL} 
\end{equation*} 
and 
\begin{equation*} 
\beta_1 :=\dfrac{1}{2i}  (W -\overline{W} ), \quad 
\beta_2 :=\dfrac{1}{2}  (Y +\overline{Y} ), \quad 
\beta_3 :=\dfrac{1}{2i} (X -\overline{X} ). 
\label{beta123tL} 
\end{equation*} 
Then $\alpha_k$, $\beta_k$ ($k=1, 2, 3$) satisfy $S_v -T_u =ST-TS$ 
for $S$, $T$ in \eqref{STtL} with $\mu_1 =\gamma_u$ and $\mu_2 =\gamma_v$. 
Therefore there exists an immersion $F$ of a neighborhood of a point 
of the $uv$-plane into $N$ with $K\equiv L_0$, $R^{\perp} \equiv 0$ 
such that the second fundamental form of $F$ does not satisfy 
the linearly dependent condition. 

Hence we obtain 

\begin{thm}\label{thm:notldL0tL} 
Let $F:M\longrightarrow N$ be a time-like and conformal immersion 
satisfying $K\equiv L_0$, $R^{\perp} \equiv 0$ and $k\not= 0, \pm i$. 
Then $F$ satisfies $X^2 \overline{Y}^2 -\overline{X}^2 Y^2 \not= 0$ 
if and only if 
$\lambda$ with $\lambda_{uu} -\lambda_{vv} +L_0 e^{2\lambda} =0$, 
$\gamma$, $k$, $\rho$, $f$ satisfy 
\eqref{fufv}, \eqref{fufvtL}, \eqref{overlinektL}, \eqref{rhotL} 
and \eqref{dgamma1_2tL}. 
In addition, 
for functions $\lambda$, $\gamma$, $k$, $\rho$, $f$ 
as above defined on a simply connected open set $O$ 
of the $uv$-plane, 
there exists a time-like and conformal immersion $F$ of $O$ into $N$ 
with $K\equiv L_0$, $R^{\perp} \equiv 0$ 
such that the second fundamental form does not satisfy 
the linearly dependent condition, 
which is unique up to an isometry of $N$. 
\end{thm} 

In the following, suppose $L_0 =0$. 
Then we can suppose $\lambda \equiv 1$. 
Then \eqref{dgamma1_2tL} is represented as 
\begin{equation} 
   \left[ 
   \begin{array}{c} 
    \gamma_u \\ 
    \gamma_v 
     \end{array} 
   \right] = 
  -\dfrac{1}{k^2 +1} 
   \left[ 
   \begin{array}{c} 
    k_u \\ 
    k_v 
     \end{array} 
   \right] 
 -i\delta 
   \Phi 
   \left[ 
   \begin{array}{c} 
    f_u \\ 
    f_v 
     \end{array} 
   \right] , \quad 
\Phi :=\dfrac{\dfrac{\partial (\overline{f}, \rho )}{\partial (u, v)}}{
              \dfrac{\partial (f, \overline{f} )}{\partial (u, v)}} . 
\label{dgamma2tL} 
\end{equation} 
Computing the first and second terms of the right side of \eqref{dgamma2tL} 
by \eqref{overlinek2+1}, \eqref{overlinekatL} and \eqref{rhotL}, 
we see that 
the imaginary part of the right side of \eqref{dgamma2tL} vanishes: 
\begin{equation} 
 {\rm Re}\,(\Phi f_a )
=-\delta\,{\rm Im}\,\dfrac{k_a}{k^2 +1} 
=-\dfrac{1}{2} \rho_a 
\label{rhoa} 
\end{equation} 
for $a=u, v$. 
Since $\gamma_{uv} =\gamma_{vu}$, we obtain $d\Phi \wedge df=0$. 
Therefore referring to discussions in Subsection~\ref{subsect:sLnpnvf}, 
we see that there exists a holomorphic function $\Tilde{\xi}$ 
of $s+it$ with \eqref{st} satisfying $\Phi =\Tilde{\xi} (f)$, i.e., 
\begin{equation}
 \dfrac{\partial (\overline{f} , \rho )}{\partial (u, v)}
=\Tilde{\xi} (f)
 \dfrac{\partial (f, \overline{f} )}{\partial (u, v)} . 
\label{gammauv2tL}
\end{equation}
By \eqref{rhoa}, we have 
\begin{equation} 
\rho =-2\,{\rm Re}\,(\Tilde{\Xi} (s+it)). 
\label{TildexiftL} 
\end{equation} 
As in Subsection~\ref{subsect:sLnpnvf}, 
\eqref{rhotL} is rewritten into the following system: 
\begin{equation} 
\begin{split} 
\delta (\sqrt{2} \varepsilon -\,{\rm cosh}\,\rho )s_u 
      -(      1              -\,{\rm sinh}\,\rho )s_v & =0, \\ 
\delta (\sqrt{2} \varepsilon +\,{\rm cosh}\,\rho )t_u 
      -(      1              +\,{\rm sinh}\,\rho )t_v & =0 
\end{split} 
\label{rhost2} 
\end{equation} 
for $\varepsilon =1$ or $-1$. 

Let $\Tilde{\Xi}$ be a holomorphic function of one complex variable. 
Then referring to \eqref{psist2}, 
we obtain a unique analytic solution $(s, t)$ 
of the initial value problem for the system \eqref{rhost2}. 
Let $f$ be a complex-valued function of two variables $u$, $v$ 
satisfying \eqref{fufv}, \eqref{fufvtL} and \eqref{rhost2} 
with \eqref{st} and \eqref{TildexiftL}. 
Then $AA' \not= 0$, and $C=f^2_u -f^2_v$ is real-valued and nonzero. 
In addition, $f$ satisfies \eqref{gammauv2tL}. 
Let $k$ be a complex-valued function 
satisfying \eqref{overlinektL}. 
Then we have $k\not= 0$, $\pm i$, and it satisfies \eqref{rhoa}. 
Therefore there exists a real-valued function $\gamma$ 
satisfying \eqref{dgamma2tL}. 
Then $\lambda \equiv 1$, $\gamma$, $k$, $\rho$, $f$ are functions 
satisfying \eqref{dgamma1_2tL}. 

Therefore from Theorem~\ref{thm:notldL0tL}, we obtain 

\begin{thm}\label{thm:notldtL} 
Let $F:M\longrightarrow E^4_1$ be a time-like and conformal immersion 
satisfying $K\equiv 0$, $R^{\perp} \equiv 0$ and $k\not= 0, \pm i$. 
Then $F$ satisfies $X^2 \overline{Y}^2 -\overline{X}^2 Y^2 \not= 0$ 
if and only if 
the complex-valued function $f$ satisfies 
\eqref{fufv}, \eqref{fufvtL} and \eqref{rhost2} 
with \eqref{st} and \eqref{TildexiftL} 
for a holomorphic function $\Tilde{\Xi}$. 
In addition, 
for functions $f$ and $\rho$ as above 
defined on a simply connected open set $O$ of the $uv$-plane, 
there exists a time-like and conformal immersion $F$ of $O$ into $E^4_1$ 
with $K\equiv 0$, $R^{\perp} \equiv 0$ 
such that the second fundamental form does not satisfy 
the linearly dependent condition. 
\end{thm} 

\section*{Acknowledgements} 

The authors are grateful to Professors Kazuyuki Enomoto, 
Atsushi Fujioka and Hideya Hashimoto for valuable comments. 
Naoya Ando is supported by 
JSPS KAKENHI Grant Number JP21K03228.

\vspace{4mm} 

\par\noindent 
\footnotesize{Naoya Ando \\ 
              Faculty of Advanced Science and Technology, 
              Kumamoto University \\ 
              2--39--1 Kurokami, Chuo-ku, Kumamoto 860--8555 Japan} 

\par\noindent  
\footnotesize{E-mail address: andonaoya@kumamoto-u.ac.jp} 

\vspace{4mm} 

\par\noindent 
\footnotesize{Ryusei Hatanaka \\ 
              Graduate School of Science and Technology, 
              Kumamoto University \\ 
              2--39--1 Kurokami, Chuo-ku, Kumamoto 860--8555 Japan} 

\par\noindent  
\footnotesize{E-mail address: 236d8010@st.kumamoto-u.ac.jp} 

\end{document}